%% file: x-1-p-2.tex
\documentclass[12pt,twoside,leqno,openany]{amsart}
\usepackage{amssymb,amsbsy,amsmath,amsfonts,amssymb,amscd,times,
graphics,color,xypic,footmisc,fancyhdr,multicol,fancybox,
graphicx,mathrsfs,rotating,ifthen,wasysym}
\usepackage[all]{xy}
\usepackage[french]{babel}
\usepackage[utf8]{inputenc}
\usepackage[T1]{fontenc}
\definecolor{black}{cmyk}{1.,1.,1.,1.0}
\definecolor{blue}{cmyk}{1.,1.,0.,0.63}
\definecolor{red}{cmyk}{0.,1.,1.,0.63}
\definecolor{green}{cmyk}{1.,0.,1.,0.63}

\input macros.tex

\let\mathcal\mathscr

\begin{document}

$\:$

\bigskip\bigskip

\begin{center}

{\large\bf Courbes projectives extrinsèques $X^1 \subset \P^2(\C)$:}

\medskip

{\large\bf harmonie avec la cohomologie intrinsèque}

\end{center}

\bigskip

\begin{center}
Jo\"el {\sc Merker}
\end{center}

\medskip

\begin{center}
\begin{minipage}[t]{11.75cm}
\baselineskip =0.35cm {\scriptsize

\centerline{\bf Table des matières}

\smallskip

{\bf 1.~Introduction
\dotfill~\pageref{introduction}.}

{\bf 2.~Géométrie initiale
\dotfill~\pageref{geometrie-initiale}.}

{\bf 3.~Fibré de jets d'ordre $\kappa \geqslant 1$ en coordonnées
\dotfill~\pageref{cotangent-en-coordonnees}.}

{\bf 4.~Construction de différentielles de jets holomorphes par élimination
\dotfill~\pageref{jets-elimination}.}

{\bf 5.~Transmission de la symétrie des différentielles de jets holomorphes
\dotfill~\pageref{jets-symetriques}.}

{\bf 6.~Annulations sur la droite à l'infini $\P_\infty^1$
\dotfill~\pageref{annulations-infini}.}

{\bf 7.~Amplitude génératrice
\dotfill~\pageref{amplitude-generatrice}.}

}\end{minipage}
\end{center}

\medskip

\section{\bf Introduction}
\label{introduction}
\HEAD{\ref{introduction}.~Introduction}{
Jo\"el Merker, Département de Mathématiques d'Orsay}

\medskip

Sur une courbe algébrique complexe projective géométriquement
lisse:
\[
X^1
\,\subset\,
\P^2(\C),
\]
pour un ordre de jets fini quelconque $\kappa \geqslant 1$,
et pour un degré homogène $m \geqslant 1$ arbitraire,  
le fibré des jets de Green-Griffiths\footnote{{\em Voir}~\cite{
Merker-2010} pour une présentation détaillée.}:
\[
\xymatrix{
E_{\kappa,m}^{\rm GG}T_{X^1}^*
\dto
\\
X^1,
}
\]
---\,\,construit fibre à fibre comme polynomialisation $m$-homogène du fibré
des jets:
\[
J^\kappa\big(\D,\,X^1\big)
\]
d'applications holomorphes locales du disque unité $\D \subset \C$ à
valeurs dans $X^1$\,\,---, est un fibré vectoriel holomorphe de rang égal à:
\[
\Card\,
\Big\{
\big(m_1,m_2,\dots,m_\kappa\big)
\in
\N^\kappa\colon\,
m_1+2\,m_2+\cdots+\kappa\,m_\kappa
=
m
\Big\}
\]
qui admet de plus une certaine filtration naturelle dont
le fibré gradué associé:
\[
{\sf Gr}^\bullet
E_{\kappa,m}^{\rm GG}T_{X^1}^*
\,\cong\,
\bigoplus_{m_1+\cdots+\kappa m_\kappa=m
\atop
m_1\geqslant 0,\,\dots,\,m_\kappa\geqslant 0}\,
\Big(
\Sym^{m_1}T_X^*
\otimes\cdots\otimes
\Sym^{m_\kappa}T_X^*
\Big)
\] 
est isomorphe à une somme
directe combinatoire
de fibrés en droites de $m_\lambda$-différentielles symétriques:
\[
\aligned
\Sym^{m_\lambda}T_X^*
&
\,\cong\,
\big(T_X^*\big)^{\otimes m_\lambda}
\\
&
\,\cong\,
\mathcal{O}_X\big(m_\lambda(d-3)\big)
\endaligned
\]
eux-mêmes isomorphes à des fibrés en droites canoniques ambiants:
\[
\mathcal{O}_X(t)
=
\mathcal{O}_{\P^2}(t)
\big\vert_X,
\]
grâce à la formule dite {\sl d'adjonction}:
\[
T_X^*
\cong
\mathcal{O}_X(d-3)
\]
dont la connaissance est anciennement établie.

\medskip

Donc:
\[
{\sf Gr}^\bullet
E_{\kappa,m}^{\rm GG}T_{X^1}^*
\,\cong\,
\bigoplus_{m_1+\cdots+\kappa m_\kappa=m
\atop
m_1\geqslant 0,\,\dots,\,m_\kappa\geqslant 0}\,
\mathcal{O}_X
\Big(
(m_1+\cdots+m_\kappa)\,(d-3)
\Big).
\]

Sachant que pour $t \geqslant d$ entier:
\[
\dim\,H^0\big(X,\mathcal{O}_X(t)\big)
\,=\,
\binom{t+2}{2}
-
\binom{t-d+2}{2},
\]
il vient en bornant
la cohomologie $H^1$, pourvu que $d \geqslant 4$ et que $m \gg 1$ soit assez
grand:
\[
\!\!\!\!\!\!\!\!\!\!\!\!\!\!\!\!\!\!\!\!\!\!\!\!
\!\!\!\!\!\!\!\!\!\!\!\!\!\!\!\!\!\!\!\!\!\!\!\!
\dim\,
H^0
\big(X,\,
E_{\kappa,m}^{\rm GG}T_X^*
\big)
\,=\,
\sum_{m_1+\cdots+\kappa m_\kappa=m}\,
\bigg\{
\binom{(m_1+\cdots+m_\kappa)(d-3)+2}{2}
-
\binom{(m_1+\cdots+m_\kappa)(d-3)-d+2}{2}
\bigg\},
\]
et un calcul simple montre que
lorsque $d \gg 1$ est grand et lorsque $m \to \infty$, 
il existe une minoration asymptotique sympathique:
\[
\dim\,
H^0
\big(X,\,
E_{\kappa,m}^{\rm GG}T_X^*
\big)
\,\geqslant\,
\frac{m^\kappa}{\kappa!\,\kappa!}\,
\Big[
d^2\,\log\,\kappa
+
d^2\,{\rm O}(1)
+
{\rm O}(d)
\Big]
+
{\rm O}\big(m^{\kappa-1}\big).
\]

\medskip

Toutefois, en décidant de passer au fibré gradué associé, cette
approche standard occulte le caractère {\em non-linéaire} essentiel
des sections holomorphes globales, elle nie aussi complètement
l'exigence mathématique ontologique incontournable de devoir embrasser
les objets géométriques modaux dans leur dépendance fondamentale aux
êtres définitionnels initiaux, et de plus, elle passe entièrement sous
silence le fait que toute la théorie des surfaces de Riemann compactes
gagne fort\,\,---\,\,en ampleur, en cohérence, en compréhension 
synthétique\,\,---\,\,à
développer une complémentarité systématique harmonieuse entre les
aspects {\em extrinsèques} et les aspects {\em intrinsèques}, comme
l'ont démontré les leçons~\cite{ Griffiths-1989} de Phillip
Griffiths, dont le résultat suivant s'inspire.

\begin{Theoreme}

\'Etant donné un ordre de jets arbitraire:
\[
\kappa
\geqslant
1,
\]
sur une courbe algébrique projective géométriquement lisse quelconque:
\[
X^1
\,\subset\,
\P^2(\C)
\]
de degré:
\[
d
\,\geqslant\,
\kappa+3,
\]
représentée\,\,---\,\,dans un système de coordonnées affines:
\[
(x,y)
\,\in\,
\C^2
\subset
\P^2
\]
associées à des coordonnées homogènes:
\[
\aligned
{}
&
[T\colon X\colon Y]
\,\in\,\P^2
\\
&
x=\frac{X}{T},
\ \ \ \ \ 
y=\frac{Y}{T},
\ \ \ \ \
\big\{T\neq 0\big\}
\cong
\C^2
\subset
\P^2,
\endaligned
\]
géométriquement adaptées pour que:
\[
\aligned
\infty_x
&
:=
[0\colon 1\colon 0]
\not\in
X^1,
\\
\infty_y
&
:=
[0\colon 0\colon 1]
\not\in
X^1,
\\
\P_\infty^1
&
:=
\big\{[0\colon X\colon Y]\big\}
\ \ 
\text{intersecte $X^1$ transversalement en $d$ points distincts},
\endaligned
\]
---\,\,comme lieu des zéros 
\[
\Big\{
(x,y)\in\C^2\colon\,
R(x,y)=0
\Big\}
\]
d'un certain polynôme $R = R(x, y)$ de degré $d \geqslant \kappa+3$
satisfaisant par lissité de $X^1 \cap \C^2$:
\[
\emptyset
=
\big\{R=0\big\}
\cap
\big\{R_x\neq 0\big\}
\cap
\big\{R_y\neq 0\big\},
\] 
alors pour tout ordre de jets intermédiaire:
\[
1
\,\leqslant\,
\lambda
\,\leqslant\,
\kappa,
\]
il existe des expressions parfaitement symétriques
à travers l'échange $x \longleftrightarrow y$:
\[
\!\!\!\!\!\!\!\!\!\!\!\!\!\!\!\!\!\!\!\!
\!\!\!\!\!\!\!\!\!\!\!\!\!\!\!\!\!\!\!\!
{\sf J}_R^\lambda
\,:=\,
\left\{
\aligned
&
\ \ \ \ \,
\frac{y^{(\lambda)}}{R_x}
+\!\!
\sum_{\mu_1+\cdots+(\lambda-1)\mu_{\lambda-1}=\lambda}
\!\!\!\!\!\!\!
\frac{\big(y'\big)^{\mu_1}\cdots
\big(y^{(\lambda-1)}\big)^{\mu_{\lambda-1}}}{
R_x}\,
\mathcal{J}_{\mu_1,\dots,\mu_{\lambda-1}}^\lambda\!
\left(
\frac{R_y}{R_x},
\bigg(
\frac{R_{x^i y^j}}{R_x}
\bigg)_{2\leqslant i+j\leqslant
\atop
\leqslant
-1+\mu_1+\cdots+\mu_{\lambda-1}}
\right),
\\
\!
&
-\,
\frac{x^{(\lambda)}}{R_y}
-\!\!
\sum_{\mu_1+\cdots+(\lambda-1)\mu_{\lambda-1}=\lambda}
\!\!\!\!\!\!\!
\frac{\big(x'\big)^{\mu_1}\cdots\big(x^{(\lambda-1)}\big)^{\mu_{\lambda-1}}}{
R_y}\,
\mathcal{J}_{\mu_1,\dots,\mu_{\lambda-1}}^\lambda\!
\left(
\frac{R_x}{R_y},
\bigg(
\frac{R_{y^ix^j}}{R_y}
\bigg)_{2\leqslant i+j\leqslant
\atop
\leqslant
-1+\mu_1+\cdots+\mu_{\lambda-1}}
\right),
\\
&
\ \ \ \ \ \ \ \ \ \ \ 
0
\ \ \ \ \ \ \ \ \ 
\text{\rm sur}\ \
X^1\cap\P_\infty^1,
\endaligned
\right.
\]
qui définissent des différentielles génératrices ${\sf J}_R^\lambda$
de $\lambda$-jets
{\em holomorphes} sur $X^1$ tout entier, notamment sur les
deux sous-ouverts de $X^1 \cap \C^2$:
\[
\aligned
&
\big\{R_x\neq 0\big\}
\ \
\text{où le fibré $J^\kappa(\D,X^1)$ est 
muni de coordonnées intrinsèques:}
\\
&
\ \ \ \ \ \ \ \ \ \ \ \ \ \ \ \ \ \ 
\big(y;y',y'',\dots,y^{(\kappa)}\big),
\\
&
\big\{R_y\neq 0\big\}
\ \
\text{où le fibré $J^\kappa(\D,X^1)$ est 
muni de coordonnées intrinsèques:}
\\
&
\ \ \ \ \ \ \ \ \ \ \ \ \ \ \ \ \ \ 
\big(x;x',x'',\dots,x^{(\kappa)}\big),
\endaligned
\]
différentielles génératrices ${\sf J}_R^\lambda$ 
qui s'annulent toutes identiquement sur
le diviseur ample $X^1 \cap \P_\infty^1$, et qui sont définies
en termes de certains polynômes:
\[
\mathcal{J}_{\mu_1,\dots,\mu_{\lambda-1}}^\lambda
=
\mathcal{J}_{\mu_1,\dots,\mu_{\lambda-1}}^\lambda
\bigg(
{\sf R}_{0,1},\,\,
\Big(
{\sf R}_{i,j}
\Big)_{2\leqslant i+j\leqslant
-1+\mu_1+\cdots+\mu_{\lambda-1}}
\bigg)
\]
à coefficients dans $\Z$ explicitables\,\,---\,\,par 
exemple pour $\lambda = 1, 2, 3$ et sur $\{ R_x \neq 0\}$:
\[
\aligned
J_R^1
&
=
\frac{y'}{R_x},
\\
J_R^2
&
=
\frac{y''}{R_x}
+
\frac{(y')^2}{R_x}
\bigg[\!
-
\frac{R_{xy}}{R_x}
+
\frac{R_y}{R_x}\,
\frac{R_{xx}}{R_x}
\bigg],
\\
J_R^3
&
=
\frac{y'''}{R_x}
+
\frac{y''y'}{R_x}
\bigg[
-\,3\,\frac{R_{xy}}{R_x}
+
3\,
\bigg(\!\frac{R_y}{R_x}\!\bigg)
\frac{R_{xx}}{R_x}
\bigg]
\,+
\\
&
\ \ \ \ \ \ \ \ \ \ \ 
+
\frac{(y')^3}{R_x}
\bigg[
-\,6\bigg(\!\frac{R_y}{R_x}\!\bigg)
\frac{R_{xy}}{R_x}\frac{R_{xx}}{R_x}
+
3\bigg(\!\frac{R_y}{R_x}\!\bigg)^2
\frac{R_{xx}}{R_x}\frac{R_{xx}}{R_x}
+
3\bigg(\!\frac{R_y}{R_x}\!\bigg)
\frac{R_{xxy}}{R_x}
-
\bigg(\!\frac{R_y}{R_x}\!\bigg)^2
\frac{R_{xxx}}{R_x}
\bigg]\text{\,\,---},
\endaligned
\]
et au moyen de ces différentielles de jets génératrices,
des sections holomorphes globales linéairement indépendantes 
du fibré des jets
de Green-Griffiths $E_{ \kappa, m}^{\rm GG} T_X^*$
sont généralement représentées comme:
\[
\boxed{\,
\sum_{m_1+2m_2+\cdots+\kappa m_\kappa=m}\,
\big({\sf J}_R^1\big)^{m_1}
\big({\sf J}_R^2\big)^{m_2}
\,\cdots\,
\big({\sf J}_R^\kappa\big)^{m_\kappa}
\cdot
{\sf G}_{m_1,m_2,\dots,m_\kappa}(x,y),\,}
\]
avec des polynômes:
\[
{\sf G}_{m_1,m_2,\dots,m_\kappa}
=
{\sf G}_{m_1,m_2,\dots,m_\kappa}(x,y)
\]
de degré:
\[
\deg\,
{\sf G}_{m_1,m_2,\dots,m_\kappa}
\,\leqslant\,
\underbrace{
m_1(d-3)
+
m_2(d-4)
+\cdots+
m_\kappa\big(d-\kappa-2\big)}_{
=:\,\delta},
\]
qui appartiennent aux espaces vectoriels quotients:
\[
\C_\delta[x,y]
\Big/
R\cdot\C_{\delta-d}[x,y],
\]
le nombre total de ces sections holomorphes linéairement indépendantes
explicites de $E_{\kappa, m}^{\rm GG} T_X^*$ étant égal à:
\[
\!\!\!\!\!\!\!\!\!\!\!\!\!\!\!\!\!\!\!\!\!\!\!\!
\!\!\!\!\!\!\!\!\!\!\!\!\!\!\!\!\!\!\!\!\!\!\!\!
\aligned
\sum_{m_1+\cdots+\kappa m_\kappa=m}
\bigg\{
\binom{m_1(d-3)+\cdots+m_\kappa(d-\kappa-2)+2}{2}
-
\binom{m_1(d-3)+\cdots+m_\kappa(d-\kappa-2)-d+2}{2}
\bigg\},
\endaligned
\]
à savoir aysmptotiquement égal à:
\[
\frac{m^\kappa}{\kappa!\,\kappa!}\,
\Big[
d^2\,\log\,\kappa
+
d^2\,{\rm O}(1)
+
{\rm O}(d)
\Big]
+
{\rm O}\big(m^{\kappa-1}\big),
\]
en harmonie et en cohérence avec la théorie intrinsèque inexplicite.\qed
\end{Theoreme}

\bigskip

\section{\bf Géométrie initiale}
\label{geometrie-initiale}
\HEAD{\ref{geometrie-initiale}.~Géométrie initiale}{
Jo\"el Merker, Département de Mathématiques d'Orsay}

\medskip

Soit $\P^2 ( \C)$ l'espace projectif complexe muni des coordonnées
homogènes:
\[
\big[T\colon X\colon Y\big].
\]
Sur les trois ouverts affines canoniques:
\[
\aligned
{\sf U}_0
&
:=
\big\{T\neq 0\big\},
\\
{\sf U}_1
&
:=
\big\{X\neq 0\big\},
\\
{\sf U}_2
&
:=
\big\{Y\neq 0\big\},
\endaligned
\]
trois systèmes de coordonnées affines canoniques:
\[
\aligned
(x_0,y_0)
:=
\big(
{\textstyle{\frac{X}{T}}},\,
{\textstyle{\frac{Y}{T}}}
\big),
\\
(x_1,y_1)
:=
\big(
{\textstyle{\frac{T}{X}}},\,
{\textstyle{\frac{Y}{X}}}
\big),
\\
(x_2,y_2)
:=
\big(
{\textstyle{\frac{X}{Y}}},\,
{\textstyle{\frac{T}{Y}}}
\big),
\endaligned
\]
sont reliés entre eux par:
\[
\aligned
x_1
&
=
\frac{1}{x_0},
\ \ \ \ \ \ \ \ \ \ \ \ \ \ \ \ \ \ \ \ \ \ \ \ \ \
y_1
=
\frac{y_0}{x_0},
\\
x_2
&
=
\frac{x_0}{y_0},
\ \ \ \ \ \ \ \ \ \ \ \ \ \ \ \ \ \ \ \ \ \ \ \ \ \
y_2
=
\frac{1}{y_0},
\\
x_2
&
=
\frac{1}{y_1},
\ \ \ \ \ \ \ \ \ \ \ \ \ \ \ \ \ \ \ \ \ \ \ \ \ \
y_2
=
\frac{x_1}{y_1}.
\endaligned
\]

Soit aussi une courbe algébrique:
\[
X^1
\,\subset\,
\P^2
\]
définie comme lieu des zéros d'un certain polynôme 
homogène:
\[
{\rm R}
=
{\rm R}
\big(
T\colon X\colon Y
\big)
\]
non identiquement nul d'un certain degré $d \geqslant 1$.

\`A ce polynôme, sont associés les trois polynômes affines
(en général non homogènes) tous de degré $\leqslant d$:
\[
\aligned
R_0(x_0,y_0)
&
:=
{\rm R}
\big(1\colon x_0\colon y_0\big),
\\
R_1(x_1,y_1)
&
:=
{\rm R}
\big(
x_1\colon 1\colon y_1\big),
\\
R_2(x_2,y_2)
&
:=
{\rm R}
\big(y_2\colon x_2\colon 1\big).
\endaligned
\]

\medskip\noindent{\bf Hypothèse géométrique.} 
La courbe $X^1 \subset \P^2$ sera toujours géométriquement lisse,
à savoir:
\[
\emptyset
=
\big\{{\rm R}=0\big\}
\cap
\big\{{\rm R}_T=0\big\}
\cap
\big\{{\rm R}_X=0\big\}
\cap
\big\{{\rm R}_Y=0\big\}.
\]
Ceci implique l'irréductibilité du polynôme homogène ${\rm R}$.

\medskip

Principalement, tout se produira dans l'une des trois cartes affines, 
par exemple celle associée à 
l'ouvert ${\sf U}_0$, c'est-à-dire en termes des coordonnées:
\[
(x_0,y_0)
\,\in\,
\C^2.
\]
Il est alors avisé d'admettre l'équivalence
notationnelle:
\[
\boxed{\,
(x,y)
\equiv
(x_0,y_0).\,}
\]

Aussi, $R = R(x, y)$ remplacera le polynôme $R_0 ( x_0, y_0)$
définissant l'équation affine de la courbe dans l'ouvert ${\sf U}_0$.

Deux points:
\[
\infty_x
\ \ \ \ \ \ \ \ \ \ \ \ \
\text{\rm et}
\ \ \ \ \ \ \ \ \ \ \ \ \
\infty_y
\]
à l'infini dans la direction de l'axe des $x$ et
à l'infini dans la direction de l'axe des $y$
existent dans $\P^2$:
\[
\infty_x
=
[0\colon 1\colon 0]
\ \ \ \ \ \ \ \ \ \ \ \ \
\text{\rm et}
\ \ \ \ \ \ \ \ \ \ \ \ \
\infty_y
=
[0\colon 0\colon 1].
\] 

L'action éventuelle d'un automorphisme
holomorphe de $\P^2$ assure (exercice) que:
\[
\infty_x
\not\in 
X^1
\ \ \ \ \ \ \ \ \ \ \ \ \
\text{\rm et}
\ \ \ \ \ \ \ \ \ \ \ \ \
\infty_y
\not\in 
X^1,
\]
et même que, après dilatation des axes de coordonnées $x$ et $y$ que:
\[
R(x,y)
=
x^d
+
a_1\,x^{d-1}y
+\cdots+
a_{d-1}\,xy^{d-1}
+
y^d
+
R_{d-1}(x,y)
+
\cdots
+
R_1(x,y)
+
R_0,
\]
où chaque polynôme $R_j(x,y)$ est homogène de degré $j$ en $(x, y)$.
{\em A posteriori} (exercice), la présence des deux monômes
$x^d$ et $y^d$ dans $R$ assure que $\infty_x \not\in X^1$
et que $\infty_y \not\in X^1$.

Comme Phillip Griffiths (\cite{ Griffiths-1989}) l'effectue régulièrement,
après action éventuelle d'un automorphisme supplémentaire de $\P^2$:
\[
\P_\infty^1
\cap 
X^1
=
\text{$d$ points distincts de multiplicité $1$}.
\]

\section{\bf Fibré de jets d'ordre $\kappa \geqslant 1$ en coordonnées}
\label{cotangent-en-coordonnees}
\HEAD{\ref{cotangent-en-coordonnees}.~Fibré de jets d'ordre 
$\kappa \geqslant 1$ en coordonnées}{
Jo\"el Merker, Département de Mathématiques d'Orsay}

Le fibré cotangent à la courbe algébrique lisse $X^1 \subset \P^2$:
\[
\xymatrix{
T_X^* \dto^\pi
\\
X,
}
\]
doit {\em impérativement} être saisi dans des trivialisations naturelles qui
sont adaptées à la disposition extrinsèque de la courbe.

Manifestement, les deux sous-ouverts naturels de l'ouvert affine ${\sf U}_0$ sont:
\[
\aligned
\big\{R_x\neq 0\big\}
&
\,\subset\,
{\sf U}_0,
\\
\big\{R_y\neq 0\big\}
&
\,\subset\,
{\sf U}_0.
\endaligned
\]

Puisque la courbe projective $\{ {\rm R} = 0 \}$ est lisse, 
quitte à effectuer encore un automorphisme de $\P^2 ( \C)$, 
le théorème de Bézout
se joue entièrement dans ${\sf U}_0$ pour les deux intersections
suivantes:
\[
\aligned
d(d-1)
&
=
\Card\,
\big\{
R=R_x=0
\big\},
\\
d(d-1)
&
=
\Card\,
\big\{
R=R_y=0
\big\},
\endaligned
\]
les points étant comptés avec multiplicité, aucune
intersection, donc, ne se faisant sur la
droite projective: 
\[
\P_\infty^1
:=
\big\{
[0\colon X\colon Y]
\big\}
\]
à l'infini.

Dans un voisinage ouvert approprié d'un point quelconque $(x_p, y_p)$
de l'ouvert semi-global $\{ R_y \neq 0\}$, le théorème
analytique des fonctions implicites représente la courbe
sous la forme d'un graphe local:
\[
\aligned
y
&
=
{\sf Y}(x)
\\
&
=
{\sf Y}_{x_p,y_p}(x),
\endaligned
\]
au moyen d'une fonction graphante analytique ${\sf Y} = {\sf Y}_{ x_p, 
y_p} ( x)$ qui dépend du point central $(x_p, y_p)$ et qui est définie
pour $\vert x - x_p \vert$ assez petit.

De même, sur l'ouvert $\{ R_x \neq 0\}$, localement au voisinage
de tout point $(x_p, y_p)$, la courbe se graphe comme:
\[
\aligned
x
&
=
{\sf X}(y)
\\
&
=
{\sf X}_{x_p,y_p}(y).
\endaligned
\]

Maintenant, si:
\[
\D
:=
\big\{
z\in\C\colon\,
\vert z\vert
<
1
\big\}
\]
désigne le disque unité ouvert de rayon
$1$ centré à l'origine dans $\C$, étant donné une application
holomorphe locale:
\[
f\colon\ \ \ 
\aligned
\D
&
\,\longrightarrow\,
\C^2
\\
\zeta
&
\,\longmapsto\,
\big(x(\zeta),y(\zeta)\big)
\endaligned
\]
alors spontanément, automatiquement\,\,---\,\,et
compulsivement\,\,---, l'application $f$ vient accompagnée de ses
différentiations formelles:
\[
\aligned
f'(\zeta)
&
=
\big(x'(\zeta),y'(\zeta)\big),
\\
f''(\zeta)
&
=
\big(x''(\zeta),y''(\zeta)\big),
\\
\cdots\cdots
&
\cdots\cdots\cdots\cdots\cdots\cdot\cdot
\\
f^{(\kappa)}(\zeta)
&
=
\big(x^{(\kappa)}(\zeta),y^{(\kappa)}(\zeta)\big),
\endaligned
\]
jusqu'à des ordres finis arbitraires $\kappa \geqslant 1$.

Dans l'ouvert affine ${\sf U}_0 \cong \C^2$, les coordonnées
de jets associées jusqu'à l'ordre 
fixé quelconque $\kappa \geqslant 1$ sont donc les
coordonnées indépendantes qui hébergent toutes ces dérivées
possibles:
\[
\Big(
x',y',\,
x'',y'',\,
\dots\dots,\,
x^{(\kappa)},y^{(\kappa)}
\Big).
\]
Ces coordonnées, intrinsèques à $\P^2$, sont extrinsèques
à la courbe $X^1 \subset \P^2$.

Mais intrinsèquement à la courbe, par exemple sur l'ouvert $\{ R_y
\neq 0\}$, puisque la coordonnée holomorphe (semi-globale) naturelle
est $x$, la coordonnées cotangente associée est $dx$\,\,---\,\,coordonnée
qu'il convient de
noter plutôt $x'$\,\,---, et généralement parlant, les coordonnées de jets
associées {\em intrinsèques} jusqu'à l'ordre $\kappa$ 
quelconque sont:
\[
\big(x',\,x'',\,\dots,\,x^{(\kappa)}\big).
\]

De même, sur l'ouvert $\{ R_x \neq 0\}$, les coordonnées
de jets associées {\em intrinsèques} jusqu'à l'ordre
$\kappa$ quelconque sont:
\[
\big(y',\,y'',\,\dots,\,y^{(\kappa)}\big).
\]

Sur l'intersection:
\[
\big\{R_y\neq 0\big\}
\cap
\big\{R_x\neq 0\big\}
\]
de ces deux ouverts, il existe, d'après la théorie générale des fibrés
holomorphes, une application de changement de trivialisation entre ces
coordonnées de jets. Comment exprimer explicitement cette
application de changement de trivialisation?

\medskip

Si donc une application holomorphe locale:
\[
f\colon\ \ \
\zeta
\,\longmapsto\,
\big(x(\zeta),y(\zeta)\big)
\]
aboutit entièrement dans la courbe:
\[
0
\equiv
R\big(x(\zeta),y(\zeta)\big),
\]
cette identité valant pour tout $\zeta \in \D$, 
une première différentiation donne alors instantanément:
\[
0
\,\equiv\,
x'(\zeta)\,
R_x\big(x(\zeta),y(\zeta)\big)
+
y'(\zeta)\,
R_y\big(x(\zeta),y(\zeta)\big).
\]

Nécessairement, l'équation ainsi obtenue doit être
interprétée en termes des coordonnées
de jets indépendantes qui hébergent les dérivées de $f$,
à savoir, comme:
\[
0
=
x'\,R_x(x,y)
+
y'\,R_y(x,y),
\]
ou, de manière équivalente, comme:
\[
y'
=
-\,
x'\,\frac{R_x}{R_y}.
\]

\begin{Lemme}
L'application de changement de trivialisation du fibré des jets
d'ordre $1$ d'applications holomorphes locales
$\D \longrightarrow X^1$ du disque unité à valeurs dans
la courbe algébrique projective $X^1 \subset \P^2$ d'équation
affine $R(x, y) = 0$ de l'ouvert
$\{ R_y \neq 0\} \times \C_{ x'}$ vers l'ouvert
$\{ R_x \neq 0\} \times \C_{ y'}$ est donnée explicitement par:
\[
\aligned
\big(
(x,y),\,\,
x'
\big)
&
\,\longmapsto\,
\big(
(x,y),\,\,
y'
\big)
\\
&\ \ \ \,
=
\bigg(
(x,y),\,\,
-\,x'\,
\frac{R_x(x,y)}{R_y(x,y)}
\bigg),
\endaligned
\]
les points-bases $(x, y) \in \C^2$ ici étant supposés appartenir
à la courbe, {\em i.e.} supposés satisfaire $R ( x, y) = 0$.\qed
\end{Lemme}

En résumé, et pour reprendre le raisonnement,
la différentiation formelle de $0 = R(x, y)$:
\[
0
=
x'\,R_x
+
y'\,R_y,
\]
s'effectue comme si 
$x$ et $y$ étaient fonction d'une
variable $\zeta \in \Delta$, et
deux résolutions équivalentes:
\[
\aligned
y'
&
=
-\,x'\,\frac{R_x}{R_y},
\\
x'
&
=
-\,y'\,
\frac{R_y}{R_x},
\endaligned
\]
expriment les changements de trivialisations, dans un sens, 
et dans l'autre sens.

\medskip

Pour connaître les changements de trivialisation entre les
deux systèmes de coordonnées de jets à l'ordre suivant $\kappa = 2$:
\[
\aligned
&
\big(x;\,x',x''\big),
\\
&
\big(y;\,y',y''\big),
\endaligned
\]
une différentiation supplémentaire:
\[
\aligned
0
&
=
x'\,R_x
+
y'\,R_y,
\\
0
&
=
x''\,R_x
+
y''\,R_y
+
(x')^2\,R_{xx}
+
2\,x'y'\,R_{xy}
+
(y')^2\,R_{yy},
\endaligned
\]
commande alors de remplacer:
\[
y'
=
-\,x'\,\frac{R_x}{R_y}
\]
dans la deuxième équation: 
\[
0
=
x''\,R_x
+
y''\,R_y
+
(x')^2\,
\bigg[
R_{xx}
-
2\,R_{xy}\,
\frac{R_x}{R_y}
+
R_{yy}\,
\bigg(
\frac{R_x}{R_y}
\bigg)^2
\bigg],
\]
et ensuite, il faut résoudre par rapport à $y''$.

\begin{Lemme}
L'application de changement de trivialisation du fibré des jets
d'ordre $2$ d'applications holomorphes locales
$\D \longrightarrow X^1$ du disque unité à valeurs dans
la courbe algébrique projective $X^1 \subset \P^2$ d'équation
affine $R(x, y) = 0$ de l'ouvert
$\{ R_y \neq 0\} \times \C_{ x', x''}^2$ vers l'ouvert
$\{ R_x \neq 0\} \times \C_{ y', y''}^2$:
\[
\aligned
\big(
(x,y),\,\,
x',\,x''
\big)
\,\longmapsto\,
\big(
(x,y),\,\,
y',\,y''
\big)
\endaligned
\]
les points-bases $(x, y) \in \C^2$ ici étant supposés appartenir
à la courbe, est donnée explicitement par:
\[
\boxed{\,
\aligned
y'
&
=
-\,x'\,\frac{R_x}{R_y},
\\
y''
&
=
-\,x''\,
\frac{R_x}{R_y}
-
(x')^2\,
\bigg[
\frac{R_{xx}}{R_x}
-
2\,
\bigg(
\frac{R_x}{R_y}
\bigg)^1\,
\frac{R_{xy}}{R_y}
+
\bigg(
\frac{R_x}{R_y}
\bigg)^2\,
\frac{R_{yy}}{R_y}
\bigg].\,\qed
\endaligned}
\]
\end{Lemme}

Ensuite, pour les jets d'ordre $3$, une troisième différentiation
est nécessaire:
\[
\aligned
0
&
=
x'''\,R_x
+
y'''\,R_y
+
\\
&
\ \ \ \ \
+
3\,x'x''\,R_{xx}
+
3\,x''y'\,R_{xy}
+
3\,x'y''\,R_{xy}
+
3\,y''y'\,R_{yy}
+
\\
&
\ \ \ \ \
+
(x')^3\,R_{xxx}
+
3\,(x')^2y'\,R_{xxy}
+
3\,x'(y')^2\,R_{xyy}
+
(y')^3\,R_{yyy}.
\endaligned
\]
La résolution par rapport à $y'''$ force à diviser par $R_y$:
\[
\aligned
y'''
&
=
-\,
x'''\,\frac{R_x}{R_y}
-
\\
&
\ \ \ \ \
-\,
3\,x''x'\,\frac{R_{xx}}{R_y}
-
3\,x''y'\,\frac{R_{xy}}{R_y}
-
3\,x'y''\,\frac{R_{xy}}{R_y}
-
3\,y''y'\,\frac{R_{yy}}{R_y}
\,-
\\
&
\ \ \ \ \
-\,
(x')^3\,\frac{R_{xxx}}{R_y}
-
3\,(x')^2y'\,\frac{R_{xxy}}{R_y}
-
3\,x'(y')^2\,\frac{R_{xyy}}{R_y}
-
(y')^3\,\frac{R_{yyy}}{R_y}.
\endaligned
\]
Mais il faut aussi remplacer les valeurs de $y'$ et de $y''$ obtenues
à l'instant:
\[
\aligned
y'''
&
=
-\,x'''\,\frac{R_x}{R_y}
\,-
\\
&
\ \ \ \ \ 
-\,
3\,x''x'\,\frac{R_{xx}}{R_y}
+
3\,x''x'\,
\bigg(\!\frac{R_x}{R_y}\!\bigg)
\frac{R_{xy}}{R_y}
+
\\
&
\ \ \ \ \ 
+
3\,x'
\bigg(
x''\,\frac{R_x}{R_y}
+
(x')^2\,
\bigg[
\frac{R_{xx}}{R_y}
-
2\,\bigg(\!\frac{R_x}{R_y}\!\bigg)
\frac{R_{xy}}{R_y}
+
\bigg(\!\frac{R_x}{R_y}\!\bigg)^2
\frac{R_{yy}}{R_y}
\bigg]
\bigg)
\frac{R_{xy}}{R_y}
\,-
\\
&
\ \ \ \ \
-\,
3\,x'\,\frac{R_x}{R_y}
\bigg(
x''\,\frac{R_x}{R_y}
+
(x')^2\,
\bigg[
\frac{R_{xx}}{R_y}
-
2\,\bigg(\!\frac{R_x}{R_y}\!\bigg)
\frac{R_{xy}}{R_y}
+
\bigg(\!\frac{R_x}{R_y}\!\bigg)^2
\frac{R_{yy}}{R_y}
\bigg]
\bigg)
\frac{R_{yy}}{R_y}
\,-
\\
&
\ \ \ \ \
-\,
(x')^3\,\frac{R_{xxx}}{R_y}
+
3\,(x')^3\,
\frac{R_x}{R_y}\,
\frac{R_{xxy}}{R_y}
-
3\,(x')^3\,
\bigg(\!\frac{R_x}{R_y}\!\bigg)^2
\frac{R_{xyy}}{R_y}
+
(x')^3
\bigg(\!\frac{R_x}{R_y}\!\bigg)^3
\frac{R_{yyy}}{R_y}.
\endaligned
\]

Simplifications et réorganisations fournissent la formule de transition
au niveau des jets d'ordre $3$:
\[
\boxed{\,
\aligned
y'''
&
=
-\,x'''\,
\frac{R_x}{R_y}
\,-
\\
&
\ \ \ \ \
-\,x''x'
\bigg[
3\,\frac{R_{xx}}{R_y}
-
6\,\bigg(\!\frac{R_x}{R_y}\!\bigg)
\frac{R_{xy}}{R_y}
+
3\bigg(\!\frac{R_x}{R_y}\!\bigg)^2
\frac{R_{yy}}{R_y}
\bigg]
\,-
\\
&
\ \ \ \ \
-\,
(x')^3
\bigg[
-\,3\,\frac{R_{xx}}{R_y}\,\frac{R_{xy}}{R_y}
+
3\bigg(\!\frac{R_x}{R_y}\!\bigg)
\frac{R_{xx}}{R_y}\,\frac{R_{yy}}{R_y}
\,+
\\
&
\ \ \ \ \ \ \ \ \ \ \ \ \ \ \ \ \ \ \ \ \,
+
6\,\bigg(\!\frac{R_x}{R_y}\!\bigg)\frac{R_{xy}}{R_y}\,
\frac{R_{xy}}{R_y}
-
9\,\bigg(\!\frac{R_x}{R_y}\!\bigg)^2
\frac{R_{xy}}{R_y}\,\frac{R_{yy}}{R_y}
+
3\bigg(\!\frac{R_x}{R_y}\!\bigg)^3
\frac{R_{yy}}{R_y}\,\frac{R_{yy}}{R_y}
\,+\,
\\
&
\ \ \ \ \ \ \ \ \ \ \ \ \ \ \ \ \ \ \ \ \,
+
\frac{R_{xxx}}{R_y}
-
3\bigg(\!\frac{R_x}{R_y}\!\bigg)
\frac{R_{xxy}}{R_y}
+
3\bigg(\!\frac{R_x}{R_y}\!\bigg)^2
\frac{R_{xyy}}{R_y}
-
\bigg(\!\frac{R_x}{R_y}\!\bigg)^3
\frac{R_{yyy}}{R_y}
\bigg].
\endaligned}
\]

\medskip

Une vision inductive autre de ces calculs les rend plus directs.

\medskip

En effet, partant de la formule de $1$-transition:
\[
y'
=
-\,x'\,
\frac{R_x}{R_y},
\]
une différentiation donne:
\[
y''
=
-\,x''\,
\frac{R_x}{R_y}
-
x'\,
\frac{x'\,R_{xx}+\boxed{y'}\,R_{xy}}{R_y}
+
\frac{x'\,R_x\big(x'\,R_{xy}+\boxed{y'}\,R_{yy}\big)}{R_y\,R_y},
\]
et il faut remplacer les $y'$ qui apparaissent:
\[
y''
=
-\,x''\,\frac{R_x}{R_y}
-
(x')^2\,\frac{R_{xx}}{R_y}
+
(x')^2\,\frac{R_x}{R_y}\,\frac{R_{xy}}{R_y}
+
(x')^2\,\frac{R_x}{R_y}\,\frac{R_{xy}}{R_y}
-
(x')^2\,\bigg(\!\frac{R_x}{R_y}\!\bigg)^2
\frac{R_{yy}}{R_y},
\]
ce qui aboutit bien à la même formule (exercice visuel).

\medskip

L'intérêt de cette approche équivalente, 
c'est qu'un seul remplacement doit être effectué 
à chaque étape.

\medskip

En effet, en partant de la formule de $2$-transition obtenue
à l'instant, et réorganisée comme il se doit:
\[
y''
=
-\,x''\,
\frac{R_x}{R_y}
-
(x')^2\,
\bigg[
\frac{R_{xx}}{R_y}
-
2\,\bigg(\!\frac{R_x}{R_y}\!\bigg)
\frac{R_{xy}}{R_y}
+
\bigg(\!\frac{R_x}{R_y}\!\bigg)^2
\frac{R_{yy}}{R_y}
\bigg],
\]
une différentiation ne fait apparaître que $y'$, mais aucun $y''$:
\[
\aligned
y'''
&
=
-\,x'''\,\frac{R_x}{R_y}
\,-
\\
&
\ \ \ \ \
\,-
\frac{x''\big(x'\,R_{xx}+\boxed{y'}\,R_{xy}\big)}{R_y}
+
\frac{x''\,R_x\,\big(x'\,R_{xy}+\boxed{y'}\,R_{yy}\big)}{R_y\,R_y}
\,-
\\
&
\ \ \ \ \
-\,2\,
x''x'\,
\bigg[
\frac{R_{xx}}{R_y}
-
2\,\bigg(\!\frac{R_x}{R_y}\!\bigg)
\frac{R_{xy}}{R_y}
+
\bigg(\!\frac{R_x}{R_y}\!\bigg)^2
\frac{R_{yy}}{R_y}
\bigg]
\,-
\\
&
\ \ \ \ \
-\,
(x')^2\,
\bigg[
\frac{x'\,R_{xxx}+\boxed{y'}\,R_{xxy}}{R_y}
-
\frac{R_{xx}\big(x'\,R_{xy}+\boxed{y'}\,R_{yy}\big)}{
R_y\,R_y}
\,-
\\
&
\ \ \ \ \ \ \ \ \ \ \ \ \ \ \ \ \ \ \ \ \
-\,
2\,
\frac{x'\,R_{xx}+\boxed{y'}\,R_{xy}}{R_y}\,
\frac{R_{xy}}{R_y}
+
2\,
\frac{R_x\big(x'\,R_{xy}+\boxed{y'}\,R_{yy}\big)}{
R_y\,R_y}\,
\frac{R_{xy}}{R_y}
\,-
\\
&
\ \ \ \ \ \ \ \ \ \ \ \ \ \ \ \ \ \ \ \ \
-\,2\,\frac{R_x}{R_y}
\,
\frac{\big(x'\,R_{xxy}+\boxed{y'}\,R_{xyy}\big)}{R_y}
+
2\,\frac{R_x}{R_y}
\frac{\big(x'\,R_{xy}+\boxed{y'}\,R_{yy}\big)}{R_y}\,
\frac{R_{xy}}{R_y}
\,+
\\
&
\ \ \ \ \ \ \ \ \ \ \ \ \ \ \ \ \ \ \ \ \
+
2\,\frac{R_x}{R_y}\,
\frac{\big(x'\,R_{xx}+\boxed{y'}\,R_{xy}\big)}{R_y}\,
\frac{R_{yy}}{R_y}
-
2\,\frac{R_x}{R_y}\,
\frac{R_x}{R_y}\,
\frac{\big(x'\,R_{xy}+\boxed{y'}\,R_{yy}\big)}{R_y}\,
\frac{R_{yy}}{R_y}
\,+
\\
&
\ \ \ \ \ \ \ \ \ \ \ \ \ \ \ \ \ \ \ \ \
+
\bigg(\!\frac{R_x}{R_y}\!\bigg)^2
\frac{\big(x'\,R_{xyy}+\boxed{y'}\,R_{yyy}\big)}{R_y}
-
\bigg(\!\frac{R_x}{R_y}\!\bigg)^2
\frac{R_{yy}}{R_y}
\big(
x'\,R_{xy}+\boxed{y'}\,R_{yy}
\big)
\bigg].
\endaligned
\]

Après remplacement (visuel) des $y'$ encadrés, 
et après simplification-réorganisation,
c'est bien la même formule de $3$-transition qui est reçue
après un effort moindre:
\[
\aligned
y'''
&
=
-\,x'''\,
\frac{R_x}{R_y}
\,-
\\
&
\ \ \ \ \
-\,x''x'
\bigg[
3\,\frac{R_{xx}}{R_y}
-
6\,\bigg(\!\frac{R_x}{R_y}\!\bigg)
\frac{R_{xy}}{R_y}
+
3\bigg(\!\frac{R_x}{R_y}\!\bigg)^2
\frac{R_{yy}}{R_y}
\bigg]
\,-
\\
&
\ \ \ \ \
-\,
(x')^3
\bigg[
-\,3\,\frac{R_{xx}}{R_y}\,\frac{R_{xy}}{R_y}
+
3\bigg(\!\frac{R_x}{R_y}\!\bigg)
\frac{R_{xx}}{R_y}\,\frac{R_{yy}}{R_y}
\,+
\\
&
\ \ \ \ \ \ \ \ \ \ \ \ \ \ \ \ \ \ \ \ \,
+
6\,\bigg(\!\frac{R_x}{R_y}\!\bigg)\frac{R_{xy}}{R_y}\,
\frac{R_{xy}}{R_y}
-
9\,\bigg(\!\frac{R_x}{R_y}\!\bigg)^2
\frac{R_{xy}}{R_y}\,\frac{R_{yy}}{R_y}
+
3\bigg(\!\frac{R_x}{R_y}\!\bigg)^3
\frac{R_{yy}}{R_y}\,\frac{R_{yy}}{R_y}
\,+\,
\\
&
\ \ \ \ \ \ \ \ \ \ \ \ \ \ \ \ \ \ \ \ \,
+
\frac{R_{xxx}}{R_y}
-
3\bigg(\!\frac{R_x}{R_y}\!\bigg)
\frac{R_{xxy}}{R_y}
+
3\bigg(\!\frac{R_x}{R_y}\!\bigg)^2
\frac{R_{xyy}}{R_y}
-
\bigg(\!\frac{R_x}{R_y}\!\bigg)^3
\frac{R_{yyy}}{R_y}
\bigg].
\endaligned
\]

\begin{Theoreme}
La différentiation d'une fonction-polynôme:
\[
\aligned
R
&
=
R(x,y)
\\
&
\equiv
R({\sf x}_1,{\sf x}_2)
\endaligned
\]
jusqu'à un ordre quelconque $\kappa \geqslant 1$ 
s'exprime explicitement comme:
\[
\aligned
R^{(\kappa)}
&
=
\sum_{e=1}^\kappa\,
\sum_{1\leqslant\lambda_1<\cdots<\lambda_e\leqslant\kappa}\,
\sum_{\mu_1\geqslant 1,\,\dots,\,\mu_e\geqslant 1}\,
\sum_{\mu_1\lambda_1+\cdots+\mu_e\lambda_e=\kappa}\,
\frac{\kappa!}{
(\lambda_1!)^{\mu_1}\,\mu_1!
\,\cdots\,
(\lambda_e!)^{\mu_e}\,\mu_e!}
\\
&
\ \ \ \ \
\sum_{j_1^1,\dots,j_{\mu_1}^1=1}^2\,
\cdots\,
\sum_{j_1^e,\dots,j_{\mu_e}^e=1}^2\,
R_{{\sf x}_{j_1^1}\cdots{\sf x}_{j_{\mu_1}^1}
\cdots\cdots
{\sf x}_{j_1^e}\cdots{\sf x}_{j_{\mu_e}^e}}
\\
&
\ \ \ \ \ \ \ \ \ \ \ \ \ \ \ \ \ \ \ \ \ \ \ \ \ \ \ \ \ \ \ \ 
\ \ \ \ \ \ \ \ \ \ \ \ 
{\sf x}_{j_1^1}^{(\lambda_1)}\cdots
{\sf x}_{j_{\mu_1}^1}^{(\lambda_1)}
\cdots\cdots
{\sf x}_{j_1^e}^{(\lambda_e)}\cdots
{\sf x}_{j_{\mu_e}^e}^{(\lambda_e)},
\endaligned
\]
en admettant l'équivalence notationnelle:
\[
(x,y)
\equiv
({\sf x}_1,{\sf x}_2).
\qed
\]
\end{Theoreme}

Pour une fonction-polynôme d'une seule variable:
\[
R
=
R({\sf x}),
\]
c'est la formule classique connue dite de Faà di Bruno:
\[
\aligned
R^{(\kappa)}
&
=
\sum_{e=1}^\kappa\,
\sum_{1\leqslant\lambda_1<\cdots<\lambda_e\leqslant\kappa}\,
\sum_{\mu_1\geqslant 1,\dots,\mu_e\geqslant 1}\,
\sum_{\mu_1\lambda_1+\cdots+\mu_e\lambda_e=\kappa}\,
\frac{\kappa!}{
(\lambda_1!)^{\mu_1}\,\mu_1!
\,\cdots\,
(\lambda_e!)^{\mu_e}\,\mu_e!}
\\
&
\ \ \ \ \ \ \ \ \ \ \ \ \ \ \ \ \ \ \ \ \ \ \ \ \ \ \ \ \ \ \ \ 
\ \ \ \ \ \ \ \ \ \ \ \ \ \ \ \ \ \ \ \ \ \ \ \ \ \ \ \ \ \ \ \ 
\ \ \ \ \ \ \ \ \ \ \ \ \
R_{{\sf x}^{\mu_1+\cdots+\mu_e}}\,
\big({\sf x}^{(\lambda_1)}\big)^{\mu_1}
\,\cdots\,
\big({\sf x}^{(\lambda_e)}\big)^{\mu_e},
\endaligned
\]
dont les premiers termes sont:
\[
\aligned
R'
&
=
{\sf x}'\,R_{\sf x},
\\
R''
&
=
{\sf x}''\,R_{\sf x}
+
({\sf x}')^2\,R_{\sf xx},
\\
R'''
&
=
{\sf x}'''\,R_{\sf x}
+
3\,{\sf x}''{\sf x}'\,R_{\sf xx}
+
({\sf x}')^3\,R_{\sf xxx},
\\
R''''
&
=
{\sf x}''''\,R_{\sf x}
+
4\,{\sf x}'''{\sf x}'\,R_{\sf xx}
+
3\,({\sf x}'')^2\,R_{\sf xx}
+
6\,{\sf x}''({\sf x}')^2\,R_{\sf xxx}
+
({\sf x}')^4\,R_{\sf xxxx},
\\
R'''''
&
=
{\sf x}'''''\,R_{\sf x}
+
5\,{\sf x}''''{\sf x}'\,R_{\sf xx}
+
10\,{\sf x}'''{\sf x}''\,R_{\sf xx}
+
10\,({\sf x}'')^2{\sf x}'\,R_{\sf xxx}
\,+
\\
&
\ \ \ \ \
+
15\,{\sf x}''({\sf x}')^2\,R_{\sf xxx}
+
10\,{\sf x}''({\sf x}')^3\,R_{\sf xxxx}
+
({\sf x}')^5\,R_{\sf xxxxx}.
\endaligned
\]

\section{\bf Formules générales de changement de trivialisation}
\label{changement-trivialisations}
\HEAD{\ref{changement-trivialisations}.~Formules générales de 
changement de trivialisation}{
Jo\"el Merker, Département de Mathématiques d'Orsay}

Maintenant, si par récurrence sur un certain entier
$\kappa \geqslant 1$, la formule de transition au niveau $\kappa$
s'écrit sous la forme:
\[
\aligned
y^{(\kappa)}
&
=
-\,x^{(\kappa)}\,
\frac{R_x}{R_y}
\,-
\\
&
-\,
\sum_{\mu_1+\cdots+(\kappa-1)\mu_{\kappa-1}=\kappa}\,
(x')^{\mu_1}
\cdots
\big(x^{(\kappa-1)}\big)^{\mu_{\kappa-1}}\,
P_{\mu_1,\dots,\mu_{\kappa-1}}^\kappa
\left(
\bigg(
\frac{R_{x^\alpha y^\beta}}{R_y}
\bigg)_{1\leqslant\alpha+\beta\leqslant\kappa
\atop
(\alpha,\beta)\neq(0,1)}
\right),
\endaligned
\]
au moyen de certains polynômes:
\[
P_{\mu_1,\dots,\mu_{\kappa-1}}^\kappa
=
P_{\mu_1,\dots,\mu_{\kappa-1}}^\kappa
\bigg(
\bigg(
{\sf R}_{\alpha,\beta}
\bigg)_{1\leqslant\alpha+\beta\leqslant\kappa
\atop
(\alpha,\beta)\neq(0,1)}
\bigg)
\]
à coefficients dans $\Z$ qui, par décision-renoncement,  
ne sont pas explicités plus avant,
alors une différentiation supplémentaire donne:
\[
\!\!\!\!\!\!\!\!\!\!\!\!\!\!\!\!\!\!\!\!
\!\!\!\!\!\!\!\!\!\!\!\!\!\!\!\!\!\!\!\!
\small
\aligned
y^{(\kappa+1)}
&
=
-\,x^{(\kappa+1)}\,
\frac{R_x}{R_y}
\,-
\\
&
\ \ \ \ \
-\,x^{(\kappa)}\,
\frac{\big(x'\,R_{xx}+\boxed{y'}\,R_{xy}\big)}{R_y}
+
x^{(\kappa)}\,
\frac{R_x}{R_y}\,
\frac{\big(x'\,R_{xy}+\boxed{y'}\,R_{yy}\big)}{R_y}
\,-
\\
&
\ \ \ \ \
-\,
\sum_{\mu_1+\cdots+(\kappa-1)\mu_{\kappa-1}=\kappa}\,
\sum_{\lambda=1}^{\kappa-1}\,
\big(x'\big)^{\mu_1}
\cdots
\mu_\lambda\,
\big(x^{(\lambda)}\big)^{\mu_\lambda-1}
x^{(\lambda+1)}
\cdots
\big(x^{(\kappa-1)}\big)^{\mu_{\kappa-1}}
\cdot
\\
&
\ \ \ \ \ \ \ \ \ \ \ \ \ \ \ \ \ \ \ \ \ \ \ \ \ \ \ \ \ \ \ \ \ \ \ \ \ 
\ \ \ \ \ \ \ \ \ \ \ \ \ \ \ \ \ \ \ \ \ \ \ \ \ \ \ \ \ \ \ \ \ \ \ \ \
\cdot
P_{\mu_1,\dots,\mu_{\kappa-1}}^\kappa
\left(
\bigg(
\frac{R_{x^\alpha y^\beta}}{R_y}
\bigg)_{1\leqslant\alpha+\beta\leqslant\kappa
\atop
(\alpha,\beta)\neq(0,1)}
\right)
\,-
\\
&
\ \ \ \ \
-\,
\sum_{\mu_1+\cdots+(\kappa-1)\mu_{\kappa-1}=\kappa}\,
\big(x'\big)^{\mu_1}
\cdots
\big(x^{\kappa-1}\big)^{\mu_{\kappa-1}}\,
\sum_{1\leqslant\alpha+\beta\leqslant\kappa
\atop
(\alpha,\beta)\neq(0,1)}\,
\\
&
\ \ \ \ \ 
\bigg[
\frac{x'\,R_{x^{\alpha+1}y^\beta}+
\boxed{y'}\,R_{x^\alpha y^{\beta+1}}}{R_y}
-
\frac{R_{x^\alpha y^\beta}}{R_y}\,
\frac{\big(x'\,R_{xy}+\boxed{y'}\,R_{yy}\big)}{R_y}
\bigg]\,
\frac{\partial P_{\mu_1,\dots,\mu_{\kappa-1}}^\kappa}{
\partial{\sf R}_{\alpha,\beta}}
\left(
\bigg(
\frac{R_{x^{\alpha_1}y^{\beta_1}}}{R_y}
\bigg)_{1\leqslant\alpha_1+\beta_1\leqslant\kappa
\atop
(\alpha_1,\beta_1)\neq(0,1)}
\right),
\endaligned
\]
et après remplacement de 
$y' = -\, x'\, \frac{ R_x}{ R_y}$: 
\[
\!\!\!\!\!\!\!\!\!\!\!\!\!\!\!\!\!\!\!\!
\!\!\!\!\!\!\!\!\!\!\!\!\!\!\!\!\!\!\!\!
\small
\aligned
y^{(\kappa+1)}
&
=
-\,x^{(\kappa+1)}\,
\frac{R_x}{R_y}
\,-
\\
&
\ \ \ \ \
-\,x^{(\kappa)}\,x'\,
\bigg[
\frac{R_{xx}}{R_y}
-
2\,
\bigg(\!\frac{R_x}{R_y}\!\bigg)\,\frac{R_{xy}}{R_y}
+
\bigg(\!\frac{R_x}{R_y}\!\bigg)^2\,\frac{R_{yy}}{R_y}
\bigg]
\,-
\\
&
\ \ \ \ \
-\,
\sum_{\mu_1+\cdots+(\kappa-1)\mu_{\kappa-1}=\kappa}\,
\sum_{\lambda=1}^{\kappa-1}\,
\big(x'\big)^{\mu_1}
\cdots
\mu_\lambda\,
\big(x^{(\lambda)}\big)^{\mu_\lambda-1}
x^{(\lambda+1)}
\cdots
\big(x^{(\kappa-1)}\big)^{\mu_{\kappa-1}}
\cdot
\\
&
\ \ \ \ \ \ \ \ \ \ \ \ \ \ \ \ \ \ \ \ \ \ \ \ \ \ \ \ \ \ \ \ \ \ \ \ \ 
\ \ \ \ \ \ \ \ \ \ \ \ \ \ \ \ \ \ \ \ \ \ \ \ \ \ \ \ \ \ \ \ \ \ \ \ \
\cdot
P_{\mu_1,\dots,\mu_{\kappa-1}}^\kappa
\left(
\bigg(
\frac{R_{x^\alpha y^\beta}}{R_y}
\bigg)_{1\leqslant\alpha+\beta\leqslant\kappa
\atop
(\alpha,\beta)\neq(0,1)}
\right)
\,-
\\
&
\ \ \ \ \ 
-\,
\sum_{\mu_1+\cdots+(\kappa-1)\mu_{\kappa-1}=\kappa}\,
\big(x'\big)^{\mu_1}
\cdots
\big(x^{(\kappa-1)}\big)^{\mu_{\kappa-1}}\,
x'\,
\sum_{1\leqslant\alpha+\beta\leqslant\kappa
\atop
(\alpha,\beta)\neq(0,1)}\,
\\
&
\ \ \ \ \
\bigg[
\frac{R_{x^{\alpha+1}y^\beta}}{R_y}
-
\frac{R_x}{R_y}\,
\frac{R_{x^\alpha y^{\beta+1}}}{R_y}
-
\frac{R_{x^\alpha y^\beta}}{R_y}\,
\frac{R_{xy}}{R_y}
+
\frac{R_{x^\alpha y^\beta}}{R_y}\,
\frac{R_x}{R_y}\,
\frac{R_{yy}}{R_y}
\bigg]\cdot
\\
&
\ \ \ \ \
\cdot
\frac{\partial P_{\mu_1,\dots,\mu_{\kappa-1}}^\kappa}{
\partial{\sf R}_{\alpha,\beta}}
\left(
\bigg(
\frac{R_{x^{\alpha_1}y^{\beta_1}}}{R_y}
\bigg)_{1\leqslant\alpha_1+\beta_1\leqslant\kappa
\atop
(\alpha_1,\beta_1)\neq(0,1)}
\right),
\endaligned
\]
c'est bien une expression du même type qui se réalise:
\[
\aligned
y^{(\kappa+1)}
&
=:
-\,
x^{(\kappa+1)}\,
\frac{R_x}{R_y}
\,-
\\
&
\ \ \ \ \
-\,
\sum_{\mu_1+\cdots+\kappa\mu_\kappa=\kappa+1}\,
\big(x'\big)^{\mu_1}
\cdots
\big(x^{(\kappa)}\big)^{\mu_\kappa}\,
P_{\mu_1,\dots,\mu_\kappa}^{\kappa+1}
\left(
\bigg(
\frac{R_{x^\alpha y^\beta}}{R_y}
\bigg)_{1\leqslant\alpha+\beta\leqslant\kappa+1
\atop
(\alpha,\beta)\neq(0,1)}
\right).
\endaligned
\]

\begin{Theoreme}
Pour tout ordre de jets $\kappa \geqslant 1$, 
les formules de changement de trivialisation:
\[
\big\{R_y\neq 0\big\}
\times
\C_{x',\dots,x^\kappa}^\kappa
\,\longrightarrow\,
\big\{R_x\neq 0\big\}
\times
\C_{y',\dots,y^\kappa}^\kappa,
\]
à savoir les composantes de l'application:
\[
\Big(
x',\dots,x^{(\lambda)},\dots,x^{(\kappa)}
\Big)
\,\longrightarrow\,
\Big(
y',\dots,y^{(\lambda)},\dots,y^{(\kappa)}
\Big)
\]
sont données pour tout ordre intermédiaire $1 \leqslant \lambda 
\leqslant \kappa$ par des formules du type:
\[
\!\!\!\!\!\!\!\!\!\!\!\!\!\!\!\!\!\!\!\!
\boxed{\,
\aligned
y^{(\lambda)}
&
=
-\,
x^{(\lambda)}\,\frac{R_x}{R_y}
\,-
\\
&
\ \ \ \ \
-\,
\sum_{\mu_1+\cdots+(\lambda-1)\mu_{\lambda-1}=\lambda
\atop
\mu_1\geqslant 0,\dots,\mu_{\lambda-1}\geqslant 0}\,
\big(x'\big)^{\mu_1}
\cdots
\big(x^{(\lambda-1)}\big)^{\mu_{\lambda-1}}\,
P_{\mu_1,\dots,\mu_{\lambda-1}}^\lambda\,
\left(
\bigg(
\frac{R_{x^\alpha y^\beta}}{R_y}
\bigg)_{1\leqslant\alpha+\beta\leqslant\lambda
\atop
(\alpha,\beta)\neq(0,1)}
\right),\,
\endaligned}
\]
au moyen de certains polynômes:
\[
P_{\mu_1,\dots,\mu_{\lambda-1}}^\lambda
=
P_{\mu_1,\dots,\mu_{\lambda-1}}^\lambda
\bigg(
\Big(
{\sf R}_{\alpha,\beta}
\Big)_{1\leqslant\alpha+\beta\leqslant\lambda
\atop
(\alpha,\beta)\neq(0,1)}
\bigg),\,
\]
à coefficients dans $\Z$.\qed
\end{Theoreme}

\begin{Question}
Expliciter complètement ces polynômes.
\end{Question}

\section{\bf Construction de différentielles de jets holomorphes 
par élimination}
\label{jets-elimination}
\HEAD{\ref{jets-elimination}.~Construction de différentielles 
de jets holomorphes 
par élimination}{
Jo\"el Merker, Département de Mathématiques d'Orsay}

Soit à nouveau l'équation affine polynomiale 
d'une courbe algébrique projective
$X^1 \subset \C^2 \subset \P^2$ géométriquement lisse:
\[
0
=
R(x,y).
\]
Une première différentiation:
\[
0
=
x'\,R_x+y'\,R_y.
\]
conduit à symétriser l'équation:
\[
\frac{y'}{R_x}
=
-\,\frac{x'}{R_y}.
\]

\begin{Lemme}
Dans l'ouvert $\big\{ R_x
\neq 0\big\}$, la courbe est un graphe: 
\[
x 
= 
{\sf X}(y)
\ \ \ \ \ \ \ \ \ \ \ \
\text{avec}\ \ 
0
\equiv
R\big({\sf X}(y),y\big),
\]
la coordonnée intrinsèque est $y$,
la coordonnée de jet est $y'$, et l'application:
\[
y
\,\longmapsto\,
\frac{y'}{
R_x\big({\sf X}(y),y\big)}
\] 
constitue une section holomorphe du fibré des jets d'ordre $1$ 
d'applications holomorphes de $\D$ à valeurs dans la courbe
$X^1 \subset \P^2$.

Dans l'ouvert $\big\{ R_y
\neq 0\big\}$, la courbe est un graphe:
\[
y 
= 
{\sf Y}(x)
\ \ \ \ \ \ \ \ \ \ \ \
\text{avec}\ \ 
0
\equiv
R\big(x,{\sf Y}(x)\big),
\]
la coordonnée intrinsèque est $x$,
la coordonnée de jet est $x'$, et l'application:
\[
x
\,\longmapsto\,
-\,\frac{x'}{
R_y\big(x,{\sf Y}(x)\big)}
\] 
constitue une section holomorphe du fibré des jets d'ordre $1$ 
d'applications holomorphes de $\D$ à valeurs dans la courbe
$X^1 \subset \P^2$.

Dans l'intersection des deux ouverts:
\[
\big\{R_x\neq 0\big\}
\cap
\big\{R_y\neq 0\big\},
\]
ces deux sections holomorphes coïncident:
\[
\frac{y'}{
R_x\big(\underbrace{{\sf X}({\sf Y}(x))}_{\equiv\,x},{\sf Y}(x)\big)}
\equiv
-\,\frac{x'}{
R_y\big(x,{\sf Y}(x)\big)},
\]
via le changement naturel de trivialisation:
\[
y'
=
-\,x'\,
\frac{R_x}{R_y}.
\qed
\]
\end{Lemme}

Par conséquent, puisque la partie affine de la courbe:
\[
X^1\cap\C^2
\,=\,
\{R_x\neq 0\big\}
\cup
\{R_y\neq 0\big\},
\] 
est lisse, l'une ou l'autre de ces deux expressions:
\[
\aligned
y
&
\,\longmapsto\,
\frac{y'}{
R_x\big({\sf X}(y),y\big)},
\\
x
&
\,\longmapsto\,
-\,\frac{x'}{
R_y\big(x,{\sf Y}(x)\big)},
\endaligned
\]
peut être prise pour définir une section holomorphe du fibré:
\[
J^1\big(\D,X^1\cap\P^2\big).
\]
Le comportement à l'infini de cette section sera étudié ultérieurement.

\medskip

Pour passer aux jets d'ordre $2$, une seconde différentiation
s'impose:
\[
\aligned
0
&
=
x'\,R_x+y'\,R_y,
\\
0
&
=
x''\,R_x+y''\,R_y
+
(x')^2\,R_{xx}
+
2\,x'y'\,R_{xy}
+
(y')^2\,R_{yy}.
\endaligned
\]
Après une résolution en $y'$ et en $y''$:
\[
\aligned
y'\,R_y
&
=
-\,x'\,R_x,
\\
y''\,R_y
&
=
-\,x''\,R_x
-
(x')^2\,R_{xx}
-
2x'y'\,R_{xy}
-
(y')^2\,R_{yy},
\endaligned
\]
une division commune par $R_x\, R_y$ amène à: 
\[
\aligned
\frac{y'}{R_x}
&
=
-\,\frac{x'}{R_y},
\\
\frac{y''}{R_x}
&
=
-\,\frac{x''}{R_y}
-
\frac{(x')^2}{R_y}\,
\frac{R_{xx}}{\boxed{R_x}}
-
2\,\frac{x'y'}{R_y}\,\frac{R_{xy}}{\boxed{R_x}}
-
\frac{(y')^2}{R_y}\,
\frac{R_{yy}}{\boxed{R_x}}.
\endaligned
\]
La première ligne ayant fourni une section holomorphe non
triviale, il est naturel de fixer de manière similaire un:

\medskip\noindent{\bf Objectif.}
{\em Dans la deuxième ligne, ne voir aucun 
$R_x$ en place dénominatoriale à droite.}

\medskip

Auparavant, bien sûr, des remplacements de $y'$ sont nécessaires:
\[
\aligned
\frac{y'}{R_x}
&
=
-\,\frac{x'}{R_y},
\\
\frac{y''}{R_x}
&
=
-\,\frac{x''}{R_y}
-
\frac{(x')^2}{R_y}\,
\bigg[
\frac{R_{xx}}{\boxed{R_x}}
-
2\,\frac{R_{xy}}{R_y}
+
\frac{R_x}{R_y}\,
\frac{R_{yy}}{R_y}
\bigg].
\endaligned
\]
Ici, deux $\frac{1}{ R_x}$ disparaissent, mais il reste encore un 
$\frac{ 1}{ R_x}$. 

\medskip\noindent{\bf Question.}
{\em Existe-t-il un moyen d'éliminer ce $\frac{ 1}{ R_x}$
rémanent à droite, sans introduire de $\frac{ 1}{ R_y}$ 
intempestif à gauche?}

\medskip

Oui. En effet, la première ligne, multipliée par un facteur
approprié:
\[
\frac{y'}{R_x}\,
x'\,\frac{R_{xx}}{R_x}
=
-\,
\frac{x'}{R_y}\,
x'\,\frac{R_{xx}}{R_x}
\]
permet de visualiser\,\,---\,\,une fois les deux équations mises 
en parallèle\,\,---\,\,que le terme intempestif (souligné ici):
\[
\aligned
\frac{y'x'}{R_x}\,
\frac{R_{xx}}{R_x}
&
=
-\,\underline{\frac{(x')^2}{R_y}\,
\frac{R_{xx}}{\boxed{R_x}}}
\\
\frac{y''}{R_x}
&
=
-\,\frac{x''}{R_y}
-
\underline{\frac{(x')^2}{R_y}\,
\bigg[
\frac{R_{xx}}{\boxed{R_x}}}
-
2\,\frac{R_{xy}}{R_y}
+
\frac{R_x}{R_y}\,
\frac{R_{yy}}{R_y}
\bigg],
\endaligned
\]
peut être annihilé par simple soustraction:
\[
\frac{y''}{R_x}
-
\frac{y'x'}{R_x}\,
\frac{R_{xx}}{R_x}
=
-\,\frac{x''}{R_y}
-
\frac{(x')^2}{R_y}\,
\bigg[
-
2\,\frac{R_{xy}}{R_y}
+
\frac{R_x}{R_y}\,
\frac{R_{yy}}{R_y}
\bigg].
\]

\`A présent, il est donc très satisfaisant de constater
qu'à gauche, seules des divisions par $R_x$ apparaissent,
tandis qu'à droite, seules des divisions par $R_y$
apparaissent.

Toutefois, un défaut demeure, puisque dans le membre de
gauche:
\[
\frac{y''}{R_x}
-
\frac{y'\,\boxed{x'}}{R_x}\,
\frac{R_{xx}}{R_x}
=
-\,\frac{x''}{R_y}
-
\frac{(x')^2}{R_y}\,
\bigg[
-
2\,\frac{R_{xy}}{R_y}
+
\frac{R_x}{R_y}\,
\frac{R_{yy}}{R_y}
\bigg],
\]
une coordonnée $x'$ apparaît, alors qu'au-dessus de l'ouvert
$\{ R_x \neq 0\}$, seules les coordonnées
de jets intrinsèques $(y, y')$ devraient
être utilisées.

\smallskip

Heureusement, ce problème peut être rapidement résolu en remplaçant:
\[
x'
=
-\,y'\,
\frac{R_y}{R_x},
\]
ce qui donne:
\[
\frac{y''}{R_x}
+
\frac{(y')^2}{R_x}\,
\frac{R_y}{R_x}\,
\frac{R_{xx}}{R_x}
=
-\,\frac{x''}{R_y}
-
\frac{(x')^2}{R_y}
\bigg[
-\,2\,\frac{R_{xy}}{R_y}
+
\frac{R_x}{R_y}\,
\frac{R_{yy}}{R_y}
\bigg],
\]
{\em sans introduire de $\frac{ 1}{ R_y}$ intempestif à gauche!}

\medskip

Donc l'énonciation du lemme précédent se généralise, à savoir dans
l'ouvert $\{ R_x \neq 0 \}$ sur lequel le fibré des jets d'ordre $2$:
\[
J^2\big(\D,X^1\cap\P^2\big),
\] 
est muni des coordonnées-fibres:
\[
\big(y',y''\big),
\]
l'application:
\[
y
\,\longmapsto\,
\frac{y''}{R_x\big({\sf X}(y),y\big)}
+
\frac{(y')^2}{R_x\big({\sf X}(y),y\big)}\,
\frac{R_y\big({\sf X}(y),y\big)}{R_x\big({\sf X}(y),y\big)}\,
\frac{R_xx\big({\sf X}(y),y\big)}{R_x\big({\sf X}(y),y\big)}\,
\]
constitue une section holomorphe, tandis que simultanément,
dans l'ouvert $\{ R_y \neq 0\}$, l'application:
\[
x\,\longmapsto\,
-\,\frac{x''}{R_y\big(x,{\sf Y}(x)\big)}
-
\frac{(x')^2}{R_y\big(x,{\sf Y}(x)\big)}\,
\bigg[
-\,2\,
\frac{R_{xy}\big(x,{\sf Y}(x)\big)}{R_y\big(x,{\sf Y}(x)\big)}
+
\frac{R_x\big(x,{\sf Y}(x)\big)}{R_y\big(x,{\sf Y}(x)\big)}\,
\frac{R_{yy}\big(x,{\sf Y}(x)\big)}{R_y\big(x,{\sf Y}(x)\big)}
\bigg],
\]
constitue une section holomorphe, {\em ces deux sections holomorphes
coïncidant point par point dans l'intersection 
$\{ R_x \neq 0 \} \cap \{ R_y \neq 0\}$ précisément
grâce à l'équation découverte à l'instant par élimination}.

\medskip

\`A nouveau, le comportement à l'infini de cette
section holomorphe définie en tout
point de $X^1 \cap \C^2$ sera étudié ultérieurement.

\medskip

En fait, un défaut formel très léger demeure encore. 
\`A travers l'échange des variables:
\[
x
\longleftrightarrow
y,
\]
la première formule:
\[
\frac{y'}{R_x}
=
-\,
\frac{x'}{R_y}
\]
était, à un signe <<\,$-$\,>> global près, {\em parfaitement
symétrique}, or tel n'est pas le cas de la formule obtenue
à l'instant:
\[
\frac{y''}{R_x}
+
\frac{(y')^2}{R_x}\,
\frac{R_y}{R_x}\,
\frac{R_{xx}}{R_x}
=
-\,\frac{x''}{R_y}
-
\frac{(x')^2}{R_y}
\bigg[
-\,\underline{2\,\frac{R_{xy}}{R_y}}
+
\frac{R_x}{R_y}\,
\frac{R_{yy}}{R_y}
\bigg].
\]

\medskip\noindent{\bf Observation.}
{\em Le terme souligné:
\[
\frac{(x')^2}{R_y}\,
\underline{2\,\frac{R_{xy}}{R_y}}
\]
peut {\em circuler} d'un côté à l'autre, simplement
en utilisant la formule du premier ordre:
\[
\frac{ y'}{ R_x} 
= 
- 
\,\frac{ x'}{ R_y},
\]
ce qui donne une équation parfaitement symétrique:
\[
\frac{(y')^2}{R_x}\,2\,
\frac{R_{xy}}{R_x}
=
\frac{(x')^2}{R_y}\,
2\,\frac{R_{xy}}{R_y},
\]
avec seulement $\frac{ 1}{ R_x}$ à gauche, et seulement 
$\frac{ 1}{ R_y}$ à droite.}

\medskip

Maintenant donc, il s'agit manifestement de diviser en deux
parties égales le terme supplémentaire à droite:
\[
\aligned
\frac{y''}{R_x}
+
\frac{(y')^2}{R_x}\,
\frac{R_y}{R_x}\,
\frac{R_{xx}}{R_x}
=
-\,\frac{x''}{R_y}
-
\frac{(x')^2}{R_y}
\bigg[
-\,\underline{({\bf 1}+1)\,\frac{R_{xy}}{R_y}}
+
\frac{R_x}{R_y}\,
\frac{R_{yy}}{R_y}
\bigg],
\endaligned
\]
et de faire passer, grâce au principe de circulation,
le {\bf 1} à gauche.

\begin{Proposition}
Une différentielle de jets holomorphe explicite et d'ordre 
précisément égal à $2$ existe sur la partie affine $X^1 \cap \P^2$
d'une courbe algébrique lisse quelconque
de degré $d \geqslant 1$ grâce à la formule:
\[
\boxed{\!
\aligned
\frac{y''}{R_x}
+
\frac{(y')^2}{R_x}
\bigg[\!
-
\frac{R_{xy}}{R_x}
+
\frac{R_y}{R_x}\,
\frac{R_{xx}}{R_x}
\bigg]
=
-\frac{x''}{R_y}
-
\frac{(x')^2}{R_y}
\bigg[\!
-\frac{R_{xy}}{R_y}
+
\frac{R_x}{R_y}\,
\frac{R_{yy}}{R_y}
\bigg],
\endaligned}
\]
formule qui est {\em parfaitement symétrique} à travers
l'échange des deux variables extrinsèques ambiantes:
\[
x
\longleftrightarrow
y,
\]
sachant qu'à gauche, seules des divisions $\frac{ 1}{ R_x}$
apparaissent, et qu'à droite, seules des divisions 
$\frac{ 1}{ R_y}$ apparaissent.
\end{Proposition}

Ce procédé de production de différentielles de jets
holomorphes par élimination et symétrisation se poursuit-il
au-delà, d'abord par exemple pour les jets d'ordre $\kappa = 3$?

\medskip

Oui, à condition de travailler plus.

\medskip

Tout d'abord, la différentiation de $0 = R(x,y)$ à l'ordre $3$
suivie d'une division par $R_x\, R_y$ et d'une résolution
en $y'''$ s'exprime comme:
\[
\aligned
\frac{y'''}{R_x}
&
=
-\,\frac{x'''}{R_y}
\,-
\\
&
\ \ \ \ \
-\,
3\,\frac{x''x'}{R_y}\,\frac{R_{xx}}{R_x}
-
3\,\frac{x''\,\boxed{y'}}{R_y}\,\frac{R_{xy}}{R_x}
-
3\,\frac{x'\,\boxed{y''}}{R_y}\,\frac{R_{xy}}{R_x}
-
3\,\frac{\boxed{y''y'}}{R_y}\,\frac{R_{yy}}{R_y}
\,-
\\
&
\ \ \ \ \
-\,
\frac{(x')^3}{R_y}\,\frac{R_{xxx}}{R_x}
-
3\,\frac{(x')^2\,\boxed{y'}}{R_y}\,\frac{R_{xxy}}{R_y}
-
3\,\frac{x'\,\boxed{(y')^2}}{R_y}\,\frac{R_{xyy}}{R_x}
-
\frac{\boxed{(y')^3}}{R_y}\,\frac{R_{yyy}}{R_x}.
\endaligned
\]
Or, seuls $x'$, $x''$, $x'''$ devraient apparaître à droite.
Il faut donc remplacer les $y'$ et $y''$ présents
en utilisant les formules déjà vues de changement de trivialisation:
\[
\aligned
y'
&
=
-\,x'\,
\frac{R_x}{R_y},
\\
y''
&
=
-\,x''\,
\frac{R_x}{R_y}
-
(x')^2\,
\bigg[
\frac{R_{xx}}{R_y}
-
2\,
\bigg(\!\frac{R_x}{R_y}\!\bigg)
\frac{R_{xy}}{R_y}
+
\bigg(\!\frac{R_x}{R_y}\!\bigg)^2
\frac{R_{yy}}{R_y}
\bigg],
\endaligned
\]
ce qui donne, sous forme d'abord brute:
\[
\aligned
\frac{y'''}{R_x}
&
=
-\,\frac{x'''}{R_y}
\,-
\\
&
\ \ \ \ \ 
-\,3\,
\frac{x''x'}{R_y}\,
\frac{R_{xx}}{\boxed{R_x}}
+
3\,\frac{x''x'}{R_y}\,\frac{R_{xy}}{R_y}
\,+
\\
&
\ \ \ \ \
+
3\,\frac{x'x''}{R_y}\,\frac{R_{xy}}{R_y}
+
3\,\frac{(x')^3}{R_y}\,
\bigg[
\frac{R_{xy}}{R_y}\,\frac{R_{xx}}{\boxed{R_x}}
-
2\,\frac{R_{xy}}{R_y}\,\frac{R_{xy}}{R_y}
+
\bigg(\!\frac{R_x}{R_y}\!\bigg)
\frac{R_{yy}}{R_y}\,
\frac{R_{xy}}{R_y}
\bigg]
\,-
\\
&
\ \ \ \ \
-\,
3\,\frac{x'x''}{R_y}\,\frac{R_{yy}}{R_y}\,
\bigg(\!\frac{R_x}{R_y}\!\bigg)
-
3\,
\frac{(x')^3}{R_y}\,
\bigg[
\frac{R_{xx}}{R_y}\,\frac{R_{yy}}{R_y}
-
2\,\bigg(\!\frac{R_x}{R_y}\!\bigg)
\frac{R_{xy}}{R_y}\,
\frac{R_{yy}}{R_y}
+
\bigg(\!\frac{R_x}{R_y}\!\bigg)^2
\frac{R_{yy}}{R_y}\,
\frac{R_{yy}}{R_y}
\bigg]
\,-
\\
&
\ \ \ \ \
-\,
\frac{(x')^3}{R_y}\,
\frac{R_{xxx}}{\boxed{R_x}}
+
3\,\frac{(x')^3}{R_y}\,
\frac{R_{xxy}}{R_y}
-
3\,\frac{(x')^3}{R_y}\,
\bigg(\!\frac{R_x}{R_y}\!\bigg)
\frac{R_{xyy}}{R_y}
+
\frac{(x')^3}{R_y}\,
\bigg(\!\frac{R_x}{R_y}\!\bigg)^2
\frac{R_{yyy}}{R_y},
\endaligned
\]
les divisions intempestives rémanentes: 
\[
\frac{*}{\boxed{R_x}}
\]
étant rendues visibles par encadrement.

Maintenant, pour éliminer le premier tel dénominateur intempestif,
il suffit de multiplier l'équation à l'ordre $2$ par le
facteur approprié:
\[
3\,x'\,\frac{R_{xx}}{R_x}\,
\left(
\frac{y''}{R_x}
=
-\,\frac{x''}{R_y}
-
\frac{(x')^2}{R_y}\,
\bigg[
\frac{R_{xx}}{R_x}
-
2\,\frac{R_{xy}}{R_y}
+
\bigg(\!\frac{R_x}{R_y}\!\bigg)
\frac{R_{yy}}{R_y}
\bigg]
\right),
\]
et de soustraire:
\[
\!\!\!\!\!\!\!\!\!\!\!\!\!\!\!\!\!\!\!\!
\!\!\!\!\!\!\!\!\!\!\!\!\!\!\!\!\!\!\!\!
\aligned
\frac{y'''}{R_x}
-
3\,\frac{y''x'}{R_x}\,
\frac{R_{xx}}{R_x}
&
=
-\,
\frac{x'''}{R_y}
+
\frac{x''x'}{R_y}\,
\bigg[
6\,\frac{R_{xy}}{R_y}
-
3\bigg(\!\frac{R_x}{R_y}\!\bigg)
\frac{R_{yy}}{R_y}
\bigg]
\,+
\\
&
\ \ \ \ \
+
\frac{(x')^3}{R_y}\,
\bigg[
3\,\frac{R_{xx}}{\boxed{R_x}}\,
\frac{R_{xx}}{\boxed{R_x}}
-
6\,\frac{R_{xy}}{R_y}\,
\frac{R_{xx}}{\boxed{R_x}}
+
3\,\frac{R_{xx}}{R_y}\,
\frac{R_{yy}}{R_y}
\,+
\\
&
\ \ \ \ \ \ \ \ \ \ \ \ \ \ \ \ \ \ \ \ \ \ \ \ \ \ \ \ \ \ \ \ \ \ \ \
\ \ \ \ 
+
3\,\frac{R_{xy}}{R_y}\,
-
6\,\frac{R_{xy}}{R_y}\,\frac{R_{xy}}{R_y}
+
3\bigg(\!\frac{R_x}{R_y}\!\bigg)
\frac{R_{yy}}{R_y}\frac{R_{xy}}{R_y}
\,-
\\
&
\ \ \ \ \ \ \ \ \ \ \ \ \ \ \ \ \ \ \ \ \ \ \ \ \ \ \ \ \ \ \ \ \ \ \ \
\ \ \ \ \ \ \ \ \ \ \ \ \ \ \ \ \ \ 
-\,
3\,\frac{R_{xx}}{R_y}\,
\frac{R_{yy}}{R_y}
+
6\,\bigg(\!\frac{R_x}{R_y}\!\bigg)
\frac{R_{xy}}{R_y}\,
\frac{R_{yy}}{R_y}
\,-
\\
&
\ \ \ \ \ \ \ \ \ \ \ \ \ \ \ \ \ \ \ \ \ \ \ \ \ \ \ \ \ \ \ \ \ \ \ \
\ \ \ \ \ \ \ \ \ \ \ \ \ \ \ \ \ \ \ \ \ \ \ \ \ \ \ \ \ \ \ \ \ \ \ \
\ \ \ \ \ \ \ \ \ \ \ \ \ \ \ \ \ \ \ \ \ \ \ \ \ \ \ \ \ \ \ \ \ \ \ \
-\,3\,
\bigg(\!\frac{R_x}{R_y}\!\bigg)
\frac{R_{yy}}{R_y}\,
\frac{R_{yy}}{R_y}
\bigg]
\,-
\\
&
\ \ \ \ \
-\,
\frac{(x')^3}{R_y}\,
\frac{R_{xxx}}{\boxed{R_x}}
+
3\,\frac{(x')^3}{R_y}\,
\frac{R_{xxy}}{R_y}
-
3\,\frac{(x')^3}{R_y}\,
\bigg(\!\frac{R_x}{R_y}\!\bigg)
\frac{R_{xyy}}{R_y}
+
\frac{(x')^3}{R_y}\,
\bigg(\!\frac{R_x}{R_y}\!\bigg)^2
\frac{R_{yyy}}{R_y},
\endaligned
\]
et enfin de collecter-réorganiser:

\[
\!\!\!\!\!\!\!\!\!\!\!\!\!\!\!\!\!\!\!\!
\aligned
\frac{y'''}{R_x}
-
3\,\frac{y''x'}{R_x}\,
\frac{R_{xx}}{R_x}
&
=
-\,
\frac{x'''}{R_y}
+
\frac{x''x'}{R_y}\,
\bigg[
6\,\frac{R_{xy}}{R_y}
-
3\bigg(\!\frac{R_x}{R_y}\!\bigg)
\frac{R_{yy}}{R_y}
\bigg]
\,+
\\
&
\ \ \ \ \
+
\frac{(x')^3}{R_y}\,
\bigg[
3\,\frac{R_{xx}}{\boxed{R_x}}\,
\frac{R_{xx}}{\boxed{R_x}}
-
3\,\frac{R_{xx}}{R_y}\,
\frac{R_{xy}}{\boxed{R_x}}
-
6\,\frac{R_{xy}}{R_y}\,
\frac{R_{xy}}{R_y}
\,+
\\
&
\ \ \ \ \ \ \ \ \ \ \ \ \ \ \ \ \ \ \ \ \ \ \ \ \ \ \ \ \ \ \ \ \ \ \ \
\ \ \ \ 
+
9\,\bigg(\!\frac{R_x}{R_y}\!\bigg)
\frac{R_{xy}}{R_y}\,\frac{R_{yy}}{R_y}
-
3\bigg(\!\frac{R_x}{R_y}\!\bigg)^2
\frac{R_{yy}}{R_y}\frac{R_{yy}}{R_y}
\,-
\\
&
\ \ \ \ \ \ \ \ \ \ \ \ \ \ \ \ \ \ \ \ \ \ \
-\,
\frac{R_{xxx}}{\boxed{R_x}}
+
3\,\frac{R_{xxy}}{R_y}
-
3\bigg(\!\frac{R_x}{R_y}\!\bigg)
\frac{R_{xyy}}{R_y}
+
\bigg(\!\frac{R_x}{R_y}\!\bigg)^2
\frac{R_{yyy}}{R_y}
\bigg].
\endaligned
\]

Pour éliminer les trois nouveaux termes
intempestifs comportant des dénominateurs
$\frac{ 1}{ R_x}$, il suffit de multiplier par trois 
facteurs appropriés la formule symétrique des jets
d'ordre $1$, éventuellement élevée au carré:
\[
\aligned
(x')^2\,\frac{R_{xxx}}{R_x}
&
\left(
\frac{y'}{R_x}
=
-\,\frac{x'}{R_y}
\right),
\\
3\,(x')^2\,\frac{R_{xx}}{R_x}\,\frac{R_{xx}}{R_x}
&
\left(
\frac{y'}{R_x}
=
-\,\frac{x'}{R_y}
\right),
\\
3\,x'\,\frac{R_{xx}}{R_x}\,\frac{R_{xy}}{1}
&
\left(
\frac{(y')^2}{R_x\,R_x}
=
\frac{(x')^2}{R_y\,R_y}
\right),
\endaligned
\]
et de soustraire ou d'additionner, ce qui donne:
\[
\aligned
&
\frac{y'''}{R_x}
-
3\,\frac{y''x'}{R_x}\,
\frac{R_{xx}}{R_x}
-
\frac{y'(x')^2}{R_x}\,
\frac{R_{xxx}}{R_x}
+
3\,\frac{(x')^2y'}{R_x}\,
\frac{R_{xx}}{R_x}\,\frac{R_{xx}}{R_x}
+
3\,\frac{x'(y')^2}{R_x}\,
\frac{R_{xx}}{R_x}\,\frac{R_{xy}}{R_x}
=
\\
&
\ \ \ \ \
=
-\,\frac{x''}{R_y}
+
\frac{x''x'}{R_y}
\bigg[
6\,\frac{R_{xy}}{R_y}
-
3\bigg(\!\frac{R_x}{R_y}\!\bigg)
\frac{R_{yy}}{R_y}
\bigg]
\,+
\\
&
\ \ \ \ \
+
\frac{(x')^3}{R_y}
\bigg[
-\,6\,\frac{R_{xy}}{R_y}\,
\frac{R_{xy}}{R_y}
+
9\,\bigg(\!\frac{R_x}{R_y}\!\bigg)
\frac{R_{xy}}{R_y}\,
\frac{R_{yy}}{R_y}
-
3\bigg(\!\frac{R_x}{R_y}\!\bigg)^2
\frac{R_{yy}}{R_y}\,
\frac{R_{yy}}{R_y}
\,+
\\
&
\ \ \ \ \ \ \ \ \ \ \ \ \ \ \ \ \ \ \ \ \
+
3\,\frac{R_{xxy}}{R_y}
-
3\bigg(\!\frac{R_x}{R_y}\!\bigg)
\frac{R_{xyy}}{R_y}
+
\bigg(\!\frac{R_x}{R_y}\!\bigg)^2
\frac{R_{yyy}}{R_y}
\bigg].
\endaligned
\]

Ainsi, puisque seules des divisions $\frac{ 1}{ R_x}$
apparaissent à gauche, et puisque seules des divisions
$\frac{ 1}{ R_y}$ apparaissent à droite, cette
équation produit une différentielle de jets
holomorphe en tout point de la partie affine
$X^1 \cap \C^2$ de la courbe projective.

Toutefois, à gauche, $x'$ apparaît quatre fois
ce qui n'est pas naturel, puisqu'au-dessus de l'ouvert
$\{ R_x = 0\}$ ne devraient apparaître
que $y'$, $y''$, $y'''$. Il suffit d'effectuer le remplacement:
\[
x'
=
-\,y'\,\frac{R_y}{R_x},
\]
ce qui fournit l'équation:
\[
\!\!\!\!\!\!\!\!\!\!\!\!\!\!\!\!\!\!\!\!
\aligned
&
\frac{y'''}{R_x}
+
\frac{y''y'}{R_x}
\bigg[
3\bigg(\!\frac{R_x}{R_y}\!\bigg)
\frac{R_{xx}}{R_x}
\bigg]
+
\frac{(y')^3}{R_x}
\bigg[
-\,3\bigg(\!\frac{R_x}{R_y}\!\bigg)
\frac{R_{xx}}{R_x}\frac{R_{xy}}{R_x}
+
3\bigg(\!\frac{R_x}{R_y}\!\bigg)^2
\frac{R_{xx}}{R_x}\frac{R_{xx}}{R_x}
-
\bigg(\!\frac{R_x}{R_y}\!\bigg)^2
\frac{R_{xxx}}{R_x}
\bigg]
=
\\
&
\ \ \ \ \
=
-\,\frac{x''}{R_y}
-
\frac{x''x'}{R_y}
\bigg[
-\,6\,\frac{R_{xy}}{R_y}
+
3\bigg(\!\frac{R_x}{R_y}\!\bigg)
\frac{R_{yy}}{R_y}
\bigg]
\,-
\\
&
\ \ \ \ \
-\,
\frac{(x')^3}{R_y}
\bigg[
6\,\frac{R_{xy}}{R_y}\,
\frac{R_{xy}}{R_y}
-
9\,\bigg(\!\frac{R_x}{R_y}\!\bigg)
\frac{R_{xy}}{R_y}\,
\frac{R_{yy}}{R_y}
+
3\bigg(\!\frac{R_x}{R_y}\!\bigg)^2
\frac{R_{yy}}{R_y}\,
\frac{R_{yy}}{R_y}
\,-
\\
&
\ \ \ \ \ \ \ \ \ \ \ \ \ \ \ \ \ \ \ \ \
-\,
3\,\frac{R_{xxy}}{R_y}
+
3\bigg(\!\frac{R_x}{R_y}\!\bigg)
\frac{R_{xyy}}{R_y}
-
\bigg(\!\frac{R_x}{R_y}\!\bigg)^2
\frac{R_{yyy}}{R_y}
\bigg].
\endaligned
\]

\`A ce stade, un léger défaut demeure encore: les deux termes à gauche
et à droite du signe <<\,$=$\,>> ne sont pas symétriques.

Manifestement, le premier travail pour
symétriser cette équation doit consister à équidistribuer
à droite et à gauche le terme:
\[
-\,6\,\frac{R_{xy}}{R_y}
\]
qui, pour l'instant, est entièrement situé à droite.

\`A cette fin, il convient de multiplier l'équation connue
par deux expression égales:
\[
\underbrace{\,\,\frac{y''}{R_x}\,\,}_{
\text{\sf multiplier par}
\atop
-\,3\,y'\,\frac{R_{xy}}{R_x}}
=
\underbrace{
-\,\frac{x''}{R_y}
-
\frac{(x')^2}{R_y}\,
\frac{R_{xx}}{R_x}
-
2\,\frac{x'y'}{R_y}\,\frac{R_{xy}}{\boxed{R_x}}
-
\frac{(y')^2}{R_y}\,
\frac{R_{yy}}{R_x}}_{
\text{\sf multiplier par}
\atop
3\,x'\,\frac{R_{xy}}{R_y}}, 
\]
pour faire circuler la moitié du $6 = 3 + 3$ en question grâce à une simple
addition:
\[
\aligned
\frac{y'''}{R_x}
&
+
\frac{y''y'}{R_x}
\bigg[
-\,3\,\frac{R_{xy}}{R_x}
+
3\bigg(\!\frac{R_y}{R_x}\!\bigg)
\frac{R_{xx}}{R_x}
\bigg]
\,+
\\
&
\ \ \ \ \
+
\frac{(y')^3}{R_x}
\bigg[
-\,3\bigg(\!\frac{R_y}{R_x}\!\bigg)
\frac{R_{xx}}{R_x}\frac{R_{xy}}{R_x}
+
3\bigg(\!\frac{R_y}{R_x}\!\bigg)^2
\frac{R_{xx}}{R_x}\frac{R_{xx}}{R_x}
-
\bigg(\!\frac{R_y}{R_x}\!\bigg)^2
\frac{R_{xxx}}{R_x}
\bigg]
=
\\
&
=
-\,\frac{x'''}{R_y}
-
\frac{x''x'}{R_y}
\bigg[
(-6+3)\,\frac{R_{xy}}{R_y}
+
3\bigg(\!\frac{R_x}{R_y}\!\bigg)
\frac{R_{yy}}{R_y}
\bigg]
\,-
\\
&
\ \ \ \ \ \ \ \ \ \ \ \ \ 
-\,
\frac{(x')^3}{R_y}
\bigg[\,
3\,\frac{R_{xx}}{\boxed{R_x}}\frac{R_{xy}}{R_y}
-
6\bigg(\!\frac{R_x}{R_y}\!\bigg)
\frac{R_{xy}}{R_y}\frac{R_{yy}}{R_y}
+
3\bigg(\!\frac{R_x}{R_y}\!\bigg)^2
\frac{R_{yy}}{R_y}\frac{R_{yy}}{R_y}
\,-
\\
&
\ \ \ \ \ \ \ \ \ \ \ \ \ \ \ \ \ \ \ \ \ \ \ \ \ \
-\,
3\,\frac{R_{xxy}}{R_y}
+
3\bigg(\!\frac{R_x}{R_y}\!\bigg)
\frac{R_{xyy}}{R_y}
-
\bigg(\!\frac{R_x}{R_y}\!\bigg)^2
\frac{R_{yyy}}{R_y}
\bigg],
\endaligned
\]
mais ce calcul introduit un $\frac{ 1}{ R_x}$ intempestif à droite.

Le lecteur a certainement deviné quels étaient les derniers
gestes de calcul qu'il s'agit d'effectuer pour obtenir (exercice)
l'expression parfaitement symétrique finale:
\[
\boxed{\,
\aligned
&
\frac{y'''}{R_x}
+
\frac{y''y'}{R_x}
\bigg[
-\,3\,\frac{R_{xy}}{R_x}
+
3\,
\bigg(\!\frac{R_y}{R_x}\!\bigg)
\frac{R_{xx}}{R_x}
\bigg]
\,+
\\
&
\ \ \ \ \ \ 
+
\frac{(y')^3}{R_x}
\bigg[
-\,6\bigg(\!\frac{R_y}{R_x}\!\bigg)
\frac{R_{xy}}{R_x}\frac{R_{xx}}{R_x}
+
3\bigg(\!\frac{R_y}{R_x}\!\bigg)^2
\frac{R_{xx}}{R_x}\frac{R_{xx}}{R_x}
\,+\,
\\
&
\ \ \ \ \ \ \ \ \ \ \ \ \ \ \ \ \ \ \ \ \ \
+
3\bigg(\!\frac{R_y}{R_x}\!\bigg)
\frac{R_{xxy}}{R_x}
-
\bigg(\!\frac{R_y}{R_x}\!\bigg)^2
\frac{R_{xxx}}{R_x}
\bigg]
=
\\
&
=
-\,\frac{x'''}{R_y}
-
\frac{x''x'}{R_y}
\bigg[
-\,3\,\frac{R_{xy}}{R_y}
+
3\bigg(\!\frac{R_x}{R_y}\!\bigg)
\frac{R_{yy}}{R_y}
\bigg]
\,-
\\
&
\ \ \ \ \ \ 
-\,
\frac{(x')^3}{R_y}
\bigg[
-\,6\bigg(\!\frac{R_x}{R_y}\!\bigg)
\frac{R_{xy}}{R_y}\frac{R_{yy}}{R_y}
+
3\bigg(\!\frac{R_x}{R_y}\!\bigg)^2
\frac{R_{yy}}{R_y}\frac{R_{yy}}{R_y}
\,+
\\
&
\ \ \ \ \ \ \ \ \ \ \ \ \ \ \ \ \ \ \ \ \ \
+\,
3\bigg(\!\frac{R_x}{R_y}\!\bigg)
\frac{R_{xyy}}{R_y}
-
\bigg(\!\frac{R_x}{R_y}\!\bigg)^2
\frac{R_{yyy}}{R_y}
\bigg].
\endaligned}
\]

\section{\bf Transmission de la symétrie des différentielles 
de jets holomorphes}
\label{jets-symetriques}
\HEAD{\ref{jets-symetriques}.~Transmission de la symétrie des différentielles 
de jets holomorphes}{
Jo\"el Merker, Département de Mathématiques d'Orsay}

Maintenant, pour les jets d'ordre supérieur $\kappa \geqslant 4$, 
il est clair que l'approche par éliminations successives
devrait toujours aboutir, mais seulement au prix de calculs
de plus en plus élaborés.

Heureusement, une autre approche plus directe fait
converger sans effort vers des équations dont
les membres de gauche et de droite sont automatiquement
symétriques par rapport à l'involution:
\[
x
\,\longleftrightarrow\,
y.
\] 

En effet par exemple, une différentiation formelle de la relation:
\[
\frac{y'}{R_x}
=
-\,\frac{x'}{R_y}
\]
donne:
\[
\frac{y''}{R_x}
-
y'\,
\frac{\big(\boxed{x'}\,R_{xx}+y'\,R_{xy}\big)}{
R_x\,R_x}
=
-\,\frac{x''}{R_y}
+
x'\,
\frac{\big(x'\,R_{xy}+\boxed{y'}\,R_{yy}\big)}{
R_y\,R_y},
\]
et deux remplacements {\em symétriques}, de $x'$ à gauche et de $y'$ 
à droite, transforment cette équation formellement en:
\[
\frac{y''}{R_x}
+
\frac{(y')^2}{R_x}
\bigg[
-\,\frac{R_{xy}}{R_x}
+
\frac{R_y}{R_x}\,\frac{R_{xx}}{R_x}
\bigg]
=
-\,\frac{x''}{R_y}
-
\frac{(x')^2}{R_y}
\bigg[
-\,\frac{R_{xy}}{R_y}
+
\frac{R_x}{R_y}\,
\frac{R_{yy}}{R_y}
\bigg],
\]
ce qui est justement l'expression parfaitement symétrique obtenue auparavant
au prix de calculs d'élimination.

\medskip

Comment justifier la véracité de ce calcul?

\medskip

C'est très simple. Dans un voisinage d'un point:
\[
(x_p,y_p)
\,\in\,
\big\{R_x\neq 0\big\}
\cap
\big\{R_y\neq 0\big\},
\]
la courbe est représentable simultanément sous forme
de deux graphes analytiques locaux:
\[
\aligned
y
&
=
{\sf Y}(x),
\\
x
&
=
{\sf X}(y),
\endaligned
\]
au moyen de deux fonctions holomorphes locales ${\sf Y} = {\sf Y} (x)$
définie près de $x_p$ et ${\sf X} = {\sf X}(y)$ définie
près de $y_p$, satisfaisant:
\[
\aligned
0
&
\equiv
R\big(x,{\sf Y}(x)\big),
\\
0
&
\equiv
R\big({\sf Y}(y),y\big).
\endaligned
\]

Dans la coordonnée intrinsèque locale $y$, l'écriture abrégée:
\[
\frac{y'}{R_x}
\]
signifie en fait rigoureusement:
\[
\frac{y'}{R_x\big({\sf X}(y),y\big)},
\]
et de même:
\[
-\,\frac{x'}{R_y}
\]
signifie rigoureusement:
\[
-\,\frac{x'}{R_y\big(x,{\sf Y}(x)\big)}.
\]

Comme cela a déjà été vu, l'égalité abrégée:
\[
\frac{y'}{R_x}
=
-\,\frac{x'}{R_y}
\]
en tenant compte de:
\[
x
=
{\sf X}(y),
\]
d'où:
\[
x'
=
y'\,{\sf X}_y(y),
\]
doit alors être lue rigoureusement comme:
\[
\aligned
\frac{y'}{R_x\big({\sf X}(y),y\big)}
&
\,\equiv\,
-\,\frac{x'}{R_y\big(x,{\sf Y}(x)\big)}
\bigg\vert_{\text{\sf après changement}
\atop
\text{\sf de coordonnée}}
\\
&
\,\equiv\,
-\,\frac{y'\,{\sf X}_y(y)}{
R_y\big({\sf X}(y),
\underbrace{{\sf Y}({\sf X}(y))}_{\equiv\,y}\big)},
\endaligned
\]
ce qui est cohérent avec la différentiation de:
\[
0
\equiv
R\big({\sf X}(y),y\big),
\]
laquelle fournit:
\[
0
=
{\sf X}_y\,R_x
+
R_y.
\]

Une autre manière\,\,---\,\,encore plus naturelle 
et directe\,\,---\,\,de voir 
pourquoi le procédé formel utilisé à l'instant est rigoureux,
c'est de rappeler que l'intérêt se porte principalement sur
des applications holomorphes locales:
\[
f\colon\ \ \
\aligned
\D
&\,\longrightarrow\,
X^1\cap\C^2
\\
\zeta
&
\,\longmapsto\,
\big(x(\zeta),y(\zeta)\big),
\endaligned
\]
à savoir satisfaisant:
\[
0
\equiv
R\big(x(\zeta),y(\zeta)\big),
\]
d'où par différentiation:
\[
0
\,\equiv\,
x'(\zeta)\,R_x\big(x(\zeta),y(\zeta)\big)
+
y'(\zeta)\,R_y\big(x(\zeta),y(\zeta)\big),
\]
puis dans l'ouvert $\{ R_y \neq 0\} \cap \{ R_x \neq 0\}$:
\[
\frac{y'(\zeta)}{R_x\big(x(\zeta),y(\zeta)\big)}
\,\equiv\,
-\,
\frac{x'(\zeta)}{R_y\big(x(\zeta),y(\zeta)\big)}.
\]

Une simple différentiation par rapport à $\zeta$ fournit alors
exactement ce qui a été obtenu à l'instant, à savoir en écrivant tous
les arguments:
\[
\aligned
&
\frac{y''(\zeta)}{
R_x\big(x(\zeta),y(\zeta)\big)}
-
y'(\zeta)\,
\frac{x'(\zeta)\,R_{xx}\big(x(\zeta),y(\zeta)\big)
+
y'(\zeta)\,R_{xy}\big(x(\zeta),y(\zeta)\big)}{
R_x\big(x(\zeta),y(\zeta)\big)\,R_x\big(x(\zeta),y(\zeta)\big)}
\,=\,
\\
&
\ \ \ \ \
\,=\,
-\,\frac{x''(\zeta)}{
R_y\big(x(\zeta),y(\zeta)\big)}
+
x'(\zeta)\,
\frac{x'(\zeta)\,R_{xy}\big(x(\zeta),y(\zeta)\big)
+
y'(\zeta)\,R_{yy}\big(x(\zeta),y(\zeta)\big)}{
R_y\big(x(\zeta),y(\zeta)\big)\,R_y\big(x(\zeta),y(\zeta)\big)}.
\endaligned
\]

\begin{Theoreme}
Pour tout ordre de jets $\kappa \geqslant 1$, il existe une
différentielle de jets holomorphe dans la partie
affine $X^1 \cap \P^2$ de la courbe qui est 
représentée par une équation parfaitement symétrique:
\[
\!\!\!\!\!\!\!\!\!\!\!\!\!\!\!\!\!\!\!\!
\aligned
&
\frac{y^{(\kappa)}}{R_x}
+
\sum_{\mu_1+\cdots+(\kappa-1)\mu_{\kappa-1}=\kappa
\atop
\mu_1\geqslant 0,\dots,\mu_{\kappa-1}\geqslant 0}
\!\!\!\!\!
\frac{\big(y'\big)^{\mu_1}\cdots\big(y^{(\kappa-1)}\big)^{\mu_{\kappa-1}}}{
R_x}\,
\mathcal{J}_{\mu_1,\dots,\mu_{\kappa-1}}^\kappa
\left(
\frac{R_y}{R_x},\,\,
\bigg(
\frac{R_{x^i y^j}}{R_x}
\bigg)_{2\leqslant i+j\leqslant
\atop
\leqslant
-1+\mu_1+\cdots+\mu_{\kappa-1}}
\right)
\,=\,
\\
&
\ \ \ \ \
\,=\,
-\,
\frac{x^{(\kappa)}}{R_y}
-
\sum_{\mu_1+\cdots+(\kappa-1)\mu_{\kappa-1}=\kappa
\atop
\mu_1\geqslant 0,\dots,\mu_{\kappa-1}\geqslant 0}
\!\!\!\!\!
\frac{\big(x'\big)^{\mu_1}\cdots\big(x^{(\kappa-1)}\big)^{\mu_{\kappa-1}}}{
R_y}\,
\mathcal{J}_{\mu_1,\dots,\mu_{\kappa-1}}^\kappa
\left(
\frac{R_x}{R_y},\,\,
\bigg(
\frac{R_{y^ix^j}}{R_y}
\bigg)_{2\leqslant i+j\leqslant
\atop
\leqslant
-1+\mu_1+\cdots+\mu_{\kappa-1}}
\right),
\endaligned
\]
en termes de certains polynômes:
\[
\mathcal{J}_{\mu_1,\dots,\mu_{\kappa-1}}^\kappa
=
\mathcal{J}_{\mu_1,\dots,\mu_{\kappa-1}}^\kappa
\bigg(
{\sf R}_{0,1},\,\,
\Big(
{\sf R}_{i,j}
\Big)_{2\leqslant i+j\leqslant
\atop
\leqslant
-1+\mu_1+\cdots+\mu_{\kappa-1}}
\bigg)
\]
à coefficients dans $\Z$.
\end{Theoreme}

\proof
\'Evidemment, la récurrence n'est pas sorcière puisque nul exigence
d'expliciteté n'est émise par cet énoncé.

Une différentiation d'ordre $1$ du membre de gauche de
cette formule symétrique admise par récurrence à l'ordre $\kappa$ donne:
\[
\!\!\!\!\!\!\!\!\!\!\!\!\!\!\!\!\!\!\!\!\!\!\!\!
\!\!\!\!\!\!\!\!\!\!\!\!\!\!\!\!\!\!\!\!\!\!\!\!
\aligned
&
\frac{y^{(\kappa+1)}}{R_x}
-
\frac{y^{(\kappa)}}{R_x}\,
\frac{\big(\boxed{x'}\,R_{xx}+y'\,R_{xy}\big)}{
R_x}
\,+
\\
&
\ \ \ \
+
\sum_{\mu_1+\cdots+(\kappa-1)\mu_{\kappa-1}=\kappa
\atop
\mu_1\geqslant 0,\dots,\mu_{\kappa-1}\geqslant 0}\,
\sum_{\lambda=1}^{\kappa-1}\,
\frac{\big(y'\big)^{\mu_1}\cdots
\mu_\lambda\,y^{(\lambda+1)}\,
\big(y^{(\lambda)}\big)^{\mu_\lambda-1}
\cdots
\big(y^{(\kappa-1)}\big)^{\mu_{\kappa-1}}}{
R_x}
\cdot
\mathcal{J}_{\mu_1,\dots,\mu_{\kappa-1}}^\kappa
\,-
\\
&
\ \ \ \ \
-\,
\sum_{\mu_1+\cdots+(\kappa-1)\mu_{\kappa-1}=\kappa
\atop
\mu_1\geqslant 0,\dots,\mu_{\kappa-1}\geqslant 0}\,
\frac{\big(y'\big)^{\mu_1}\cdots
\big(y^{(\kappa-1)}\big)^{\mu_{\kappa-1}}}{R_x}
\cdot
\frac{\big(\boxed{x'}\,R_{xx}+y'\,R_{xy}\big)}{R_x}
\cdot
\mathcal{J}_{\mu_1,\dots,\mu_{\kappa-1}}^\kappa
\,+
\\
&
\ \ \ \ \
+\!\!\!
\sum_{\mu_1+\cdots+(\kappa-1)\mu_{\kappa-1}=\kappa
\atop
\mu_1\geqslant 0,\dots,\mu_{\kappa-1}\geqslant 0}
\!\!\!\!\!
\frac{\big(y'\big)^{\mu_1}\cdots
\big(y^{(\kappa-1)}\big)^{\mu_{\kappa-1}}}{R_x}
\cdot
\frac{\partial\mathcal{J}_{\mu_1,\dots,\mu_{\kappa-1}}^\kappa}{
\partial{\sf R}_{0,1}}
\cdot
\bigg[
\frac{\boxed{x'}\,R_{xy}+y'\,R_{yy}}{R_x}
-
\frac{R_y}{R_x}\,
\frac{\big(\boxed{x'}\,R_{xx}+y'\,R_{xy}\big)}{R_x}
\bigg]
\,+
\\
&
\ \ \ \ \
+
\sum_{\mu_1+\cdots+(\kappa-1)\mu_{\kappa-1}=\kappa
\atop
\mu_1\geqslant 0,\dots,\mu_{\kappa-1}\geqslant 0}\,
\frac{\big(y'\big)^{\mu_1}\cdots
\big(y^{(\kappa-1)}\big)^{\mu_{\kappa-1}}}{R_x}
\cdot
\sum_{2\leqslant i+j\leqslant -1+\mu_1+\cdots+\mu_{\kappa-1}}\,
\frac{\partial\mathcal{J}_{\mu_1,\dots,\mu_{\kappa-1}}^\kappa}{
\partial{\sf R}_{i,j}}
\cdot
\\
&
\ \ \ \ \ \ \ \ \ \ \ \ \ \ \ \ \ \ \ \ \ \ \ \ \ \ \ \ \ \ \ \ \    
\cdot
\bigg[
\frac{\boxed{x'}\,R_{x^{i+1}y^j}+y'\,R_{x^iy^{j+1}}}{R_x}
-
\frac{R_{x^iy^j}}{R_x}\,
\frac{\big(\boxed{x'}\,R_{xx}+y'\,R_{xy}\big)}{R_x}
\bigg],
\endaligned
\]
et la différentiation du membre de droite donne une
expression entièrement symétrique.

Il faut alors remplacer tous les $x'$ encadrés:
\[
\!\!\!\!\!\!\!\!\!\!\!\!\!\!\!\!\!\!\!\!
\!\!\!\!\!\!\!\!\!\!\!\!\!\!\!\!\!\!\!\!
\aligned
&
\,=\,
\frac{y^{(\kappa+1)}}{R_x}
+
\frac{y^{(\kappa)}y'}{R_x}
\bigg[
\bigg(\!\frac{R_y}{R_x}\!\bigg)
\frac{R_{xx}}{R_x}
-
\frac{R_{xy}}{R_x}
\bigg]
\,+
\\
&
\ \ \ \
+
\sum_{\mu_1+\cdots+(\kappa-1)\mu_{\kappa-1}=\kappa
\atop
\mu_1\geqslant 0,\dots,\mu_{\kappa-1}\geqslant 0}\,
\sum_{\lambda=1}^{\kappa-1}\,
\mu_\lambda\,
\frac{\big(y'\big)^{\mu_1}\cdots
\big(y^{(\lambda)}\big)^{\mu_\lambda-1}
\big(y^{(\lambda+1)}\big)^{\mu_{\lambda+1}+1}
\cdots
\big(y^{(\kappa-1)}\big)^{\mu_{\kappa-1}}}{
R_x}
\cdot
\mathcal{J}_{\mu_1,\dots,\mu_{\kappa-1}}^\kappa
\,-
\\
&
-\,
\sum_{\mu_1+\cdots+(\kappa-1)\mu_{\kappa-1}=\kappa
\atop
\mu_1\geqslant 0,\dots,\mu_{\kappa-1}\geqslant 0}\,
\frac{y'\,\big(y'\big)^{\mu_1}\cdots
\big(y^{(\kappa-1)}\big)^{\mu_{\kappa-1}}}{R_x}
\cdot
\bigg[
-\,
\bigg(\!\frac{R_y}{R_x}\!\bigg)
\frac{R_{xx}}{R_x}
+
\frac{R_{xy}}{R_x}
\bigg]
\cdot
\mathcal{J}_{\mu_1,\dots,\mu_{\kappa-1}}^\kappa
\,+
\\
&
\ \ \ \ \
+
\sum_{\mu_1+\cdots+(\kappa-1)\mu_{\kappa-1}=\kappa
\atop
\mu_1\geqslant 0,\dots,\mu_{\kappa-1}\geqslant 0}\,
\frac{y'\,\big(y'\big)^{\mu_1}\cdots
\big(y^{(\kappa-1)}\big)^{\mu_{\kappa-1}}}{R_x}
\cdot
\frac{\partial\mathcal{J}_{\mu_1,\dots,\mu_{\kappa-1}}^\kappa}{
\partial{\sf R}_{0,1}}
\cdot
\bigg[
-\,2
\bigg(\!\frac{R_y}{R_x}\!\bigg)
\frac{R_{xy}}{R_x}
+
\frac{R_{yy}}{R_x}
+
\bigg(\!\frac{R_y}{R_x}\!\bigg)^2
\frac{R_{xx}}{R_x}
\bigg]
\,+
\\
&
\ \ \ \ \
+
\sum_{\mu_1+\cdots+(\kappa-1)\mu_{\kappa-1}=\kappa
\atop
\mu_1\geqslant 0,\dots,\mu_{\kappa-1}\geqslant 0}\,
\frac{y'\,\big(y'\big)^{\mu_1}\cdots
\big(y^{(\kappa-1)}\big)^{\mu_{\kappa-1}}}{R_x}
\sum_{2\leqslant i+j\leqslant -1+\mu_1+\cdots+\mu_{\kappa-1}}\,
\frac{\partial\mathcal{J}_{\mu_1,\dots,\mu_{\kappa-1}}^\kappa}{
\partial{\sf R}_{i,j}}
\cdot
\\
&
\ \ \ \ \ \ \ \ \ \ \ \ \ \ \ \ \ \ \ \ \ \ \ \ \ \ \ \ \ \ \ \ \  
\cdot
\bigg[
-\,
\bigg(\!\frac{R_y}{R_x}\!\bigg)
\frac{R_{x^{i+1}y^j}}{R_x}
+
\frac{R_{x^iy^{j+1}}}{R_x}
+
\bigg(\!\frac{R_y}{R_x}\!\bigg)
\frac{R_{xx}}{R_x}
\frac{R_{x^iy^j}}{R_x}
-
\frac{R_{xy}}{R_x}
\frac{R_{x^iy^j}}{R_x}
\bigg],
\endaligned
\]
et effectuer aussi les remplacements symétriques de tous les
$y'$ qui apparaissent dans la différentiation
d'ordre $1$ du membre de droite. Le résultat obtenu
est bien de la forme générale annoncée, et ce, au niveau $\kappa + 1$
(exercice visuel et mental).
\endproof

Maintenant, pour éclaircir encore plus avant ce procédé 
et s'assurer de sa cohérence symbolique avec ce qui précède, il est avisé
d'examiner le résultat explicite, par exemple à l'ordre
$\kappa = 3$. 

Une différentiation de:
\[
\frac{y''}{R_x}
-
y'\,
\frac{\big(x'\,R_{xx}+y'\,R_{xy}\big)}{R_x\,R_x}
\]
donne, tous remplacements et toutes simplifications étant faits:
\[
\aligned
&
\frac{y'''}{R_x}
+
\frac{y''y'}{R_x}
\bigg[
-\,3\,
\frac{R_{xy}}{R_x}
+
3\bigg(\!\frac{R_y}{R_x}\!\bigg)
\frac{R_{xx}}{R_x}
\bigg]
\,+
\\
&
\ \ \ \ \ \
+
\frac{(y')^3}{R_x}
\bigg[
\underline{2\,\frac{R_{xy}}{R_x}
\frac{R_{xy}}{R_x}}
+
\underline{\frac{R_{xx}}{R_x}
\frac{R_{yy}}{R_x}}
-
6\,\bigg(\!\frac{R_y}{R_x}\!\bigg)
\frac{R_{xx}}{R_x}
\frac{R_{xy}}{R_x}
+
3\bigg(\!\frac{R_y}{R_x}\!\bigg)^2
\frac{R_{xx}}{R_x}
\frac{R_{xx}}{R_x}
\,-
\\
&
\ \ \ \ \ \ \ \ \ \ \ \ \ \ \ \ \ \ \ \ \ 
-\,
\frac{R_{xyy}}{R_x}
+
2\,\bigg(\!\frac{R_y}{R_x}\!\bigg)
\frac{R_{xxy}}{R_x}
-
\bigg(\!\frac{R_y}{R_x}\!\bigg)^2
\frac{R_{xxx}}{R_x}
\bigg].
\endaligned
\]
Ce n'est pas la même expression que celle obtenue auparavant
par élimination! En particulier, les termes soulignés
diffèrent! D'autes encore diffèrent!

Il y a une explication à cela: {\em ces termes soulignés
peuvent circuler d'un côté et de l'autre
du signe <<\,$=$\,>>}, sans compromettre le fait
que seules des divisions par $R_x$ sont autorisées
à droite, et seules des divisions par $R_y$ sont
autorisées à gauche: 
\[
\aligned
&
\frac{y'''}{R_x}
+
\frac{y''y'}{R_x}
\bigg[
-\,3\,
\frac{R_{xy}}{R_x}
+
3\bigg(\!\frac{R_y}{R_x}\!\bigg)
\frac{R_{xx}}{R_x}
\bigg]
\,+
\\
&
\ \ \ \ \ \
+
\frac{(y')^3}{R_x}
\bigg[
\zero{2\,\frac{R_{xy}}{R_x}
\frac{R_{xy}}{R_x}}
+
\zerozero{\frac{R_{xx}}{R_x}
\frac{R_{yy}}{R_x}}
-
6\,\bigg(\!\frac{R_y}{R_x}\!\bigg)
\frac{R_{xx}}{R_x}
\frac{R_{xy}}{R_x}
+
3\bigg(\!\frac{R_y}{R_x}\!\bigg)^2
\frac{R_{xx}}{R_x}
\frac{R_{xx}}{R_x}
\,-
\\
&
\ \ \ \ \ \ \ \ \ \ \ \ \ \ \ \ \ \ \ \ \ 
-\,
\frac{R_{xyy}}{R_x}
+
2\,\bigg(\!\frac{R_y}{R_x}\!\bigg)
\frac{R_{xxy}}{R_x}
-
\bigg(\!\frac{R_y}{R_x}\!\bigg)^2
\frac{R_{xxx}}{R_x}
\bigg]
\,=\,
\\
&
\,=\,
-\,
\frac{x'''}{R_y}
-
\frac{x''x'}{R_y}
\bigg[
-\,3\,
\frac{R_{xy}}{R_y}
+
3\bigg(\!\frac{R_x}{R_y}\!\bigg)
\frac{R_{yy}}{R_y}
\bigg]
\,+
\\
&
\ \ \ \ \ \
+
\frac{(x')^3}{R_y}
\bigg[
\zero{2\,\frac{R_{xy}}{R_y}
\frac{R_{xy}}{R_y}}
+
\zerozero{\frac{R_{yy}}{R_y}
\frac{R_{xx}}{R_y}}
-
6\,\bigg(\!\frac{R_x}{R_y}\!\bigg)
\frac{R_{yy}}{R_y}
\frac{R_{xy}}{R_y}
+
3\bigg(\!\frac{R_x}{R_y}\!\bigg)^2
\frac{R_{yy}}{R_y}
\frac{R_{yy}}{R_y}
\,-
\\
&
\ \ \ \ \ \ \ \ \ \ \ \ \ \ \ \ \ \ \ \ \ 
-\,
\frac{R_{xxy}}{R_y}
+
2\,\bigg(\!\frac{R_x}{R_y}\!\bigg)
\frac{R_{xyy}}{R_y}
-
\bigg(\!\frac{R_x}{R_y}\!\bigg)^2
\frac{R_{yyy}}{R_y}
\bigg],
\endaligned
\]
et qui plus est, les termes soulignés s'annihilent
de part et d'autre, en tenant compte de la relation $\frac{ y'}{ R_x} = 
-\, \frac{ x'}{ R_y}$, convenablement multipliée.
Après cette simplification:
\[
\aligned
&
\frac{y'''}{R_x}
+
\frac{y''y'}{R_x}
\bigg[
-\,3\,
\frac{R_{xy}}{R_x}
+
3\bigg(\!\frac{R_y}{R_x}\!\bigg)
\frac{R_{xx}}{R_x}
\bigg]
\,+
\\
&
\ \ \ \ \ \
+
\frac{(y')^3}{R_x}
\bigg[
-\,
6\,\bigg(\!\frac{R_y}{R_x}\!\bigg)
\frac{R_{xx}}{R_x}
\frac{R_{xy}}{R_x}
+
3\bigg(\!\frac{R_y}{R_x}\!\bigg)^2
\frac{R_{xx}}{R_x}
\frac{R_{xx}}{R_x}
\,-
\\
&
\ \ \ \ \ \ \ \ \ \ \ \ \ \ \ \ \ \ \ \ \ 
-\,
\frac{R_{xyy}}{R_x}
+
2\,\bigg(\!\frac{R_y}{R_x}\!\bigg)
\frac{R_{xxy}}{R_x}
-
\bigg(\!\frac{R_y}{R_x}\!\bigg)^2
\frac{R_{xxx}}{R_x}
\bigg]
\,=\,
\\
&
\,=\,
-\,
\frac{x'''}{R_y}
-
\frac{x''x'}{R_y}
\bigg[
-\,3\,
\frac{R_{xy}}{R_y}
+
3\bigg(\!\frac{R_x}{R_y}\!\bigg)
\frac{R_{yy}}{R_y}
\bigg]
\,+
\\
&
\ \ \ \ \ \
+
\frac{(x')^3}{R_y}
\bigg[
-\,
6\,\bigg(\!\frac{R_x}{R_y}\!\bigg)
\frac{R_{yy}}{R_y}
\frac{R_{xy}}{R_y}
+
3\bigg(\!\frac{R_x}{R_y}\!\bigg)^2
\frac{R_{yy}}{R_y}
\frac{R_{yy}}{R_y}
\,-
\\
&
\ \ \ \ \ \ \ \ \ \ \ \ \ \ \ \ \ \ \ \ \ 
-\,
\frac{R_{xxy}}{R_y}
+
2\,\bigg(\!\frac{R_x}{R_y}\!\bigg)
\frac{R_{xyy}}{R_y}
-
\bigg(\!\frac{R_x}{R_y}\!\bigg)^2
\frac{R_{yyy}}{R_y}
\bigg],
\endaligned
\]
et il est plus avantageux de repartir de cette identité pour différentier
une fois de plus et obtenir une différentielle de jets
d'ordre $\kappa = 4$ qui est holomorphe sur $X^1 \cap \C^2$.
 
Tous remplacements effectués et toutes simplifications faites, la
différentiation du membre de droite devient:
\[
\!\!\!\!\!\!\!\!\!\!\!\!\!\!\!\!\!\!\!\!\!\!\!\!\!\!
\aligned
&
\frac{y''''}{R_x}
+
\frac{y'''y'}{R_x}
\bigg[
-\,4\,\frac{R_{xy}}{R_x}
+
4\bigg(\!\frac{R_y}{R_x}\!\bigg)
\frac{R_{xx}}{R_x}
\bigg]
+
\frac{y''y''}{R_x}
\bigg[
-\,3\,\frac{R_{xy}}{R_x}
+
3\bigg(\!\frac{R_y}{R_x}\!\bigg)
\frac{R_{xx}}{R_x}
\bigg]
+
\\
&
+
\frac{y''(y')^2}{R_x}
\bigg[
3\,\frac{R_{xy}}{R_x}\frac{R_{xy}}{R_x}
+
3\,\frac{R_{yy}}{R_x}\frac{R_{xx}}{R_x}
-
30\bigg(\!\frac{R_y}{R_x}\!\bigg)
\frac{R_{xx}}{R_x}\frac{R_{xy}}{R_x}
+
15\,\bigg(\!\frac{R_y}{R_x}\!\bigg)^2
\frac{R_{xx}}{R_x}\frac{R_{xx}}{R_x}
\,-
\\
&
\ \ \ \ \ \ \ \ \ \ \ \ \ \ \ \
-\,
6\,\frac{R_{xyy}}{R_x}
+
12\bigg(\!\frac{R_y}{R_x}\!\bigg)
\frac{R_{xxy}}{R_x}
-
6\bigg(\!\frac{R_y}{R_x}\!\bigg)^2
\frac{R_{xxx}}{R_x}
\bigg]
\,+
\endaligned
\]
\[
\!\!\!\!\!\!\!\!\!\!\!\!\!\!\!\!\!\!\!\!\!\!\!\!\!\!
\aligned
&
+
\frac{(y')^4}{R_x}
\bigg[
-\,6\,\frac{R_{xx}}{R_x}\frac{R_{xy}}{R_x}\frac{R_{yy}}{R_x}
+
30\bigg(\!\frac{R_y}{R_x}\!\bigg)
\frac{R_{xx}}{R_x}\frac{R_{xy}}{R_x}\frac{R_{xy}}{R_x}
+
6\bigg(\!\frac{R_y}{R_x}\!\bigg)
\frac{R_{xx}}{R_x}\frac{R_{xx}}{R_x}\frac{R_{yy}}{R_x}
\,-
\\
&
\ \ \ \ \ \ \ \ \ \ \ \ \ \ \ \
-\,
45\bigg(\!\frac{R_y}{R_x}\!\bigg)^2
\frac{R_{xx}}{R_x}\frac{R_{xx}}{R_x}\frac{R_{xy}}{R_x}
+
15\bigg(\!\frac{R_y}{R_x}\!\bigg)^2
\frac{R_{xx}}{R_x}\frac{R_{xx}}{R_x}\frac{R_{xx}}{R_x}
\,+
\\
&
\ \ \ \ \ \ \ \ \ \ \ \ \ \ \ \
+
2\,\frac{R_{xy}}{R_x}\frac{R_{xyy}}{R_x}
+
2\,\frac{R_{yy}}{R_x}\frac{R_{xxy}}{R_x}
-
8\bigg(\!\frac{R_y}{R_x}\!\bigg)
\frac{R_{xx}}{R_x}\frac{R_{xyy}}{R_x}
-
14\bigg(\!\frac{R_y}{R_x}\!\bigg)
\frac{R_{xy}}{R_x}\frac{R_{xxy}}{R_x}
-
2\bigg(\!\frac{R_y}{R_x}\!\bigg)
\frac{R_{yy}}{R_x}
\frac{R_{xxx}}{R_x}
\,+
\\
&
\ \ \ \ \ \ \ \ \ \ \ \ \ \ \ \
+
18\bigg(\!\frac{R_y}{R_x}\!\bigg)^2\frac{R_{xx}}{R_x}
\frac{R_{xxy}}{R_x}
+
12\bigg(\!\frac{R_y}{R_x}\!\bigg)^2
\frac{R_{xy}}{R_x}\frac{R_{xxx}}{R_x}
-
10\bigg(\!\frac{R_y}{R_x}\!\bigg)^3
\frac{R_{xx}}{R_x}\frac{R_{xxx}}{R_x}
\,-
\\
&
\ \ \ \ \ \ \ \ \ \ \ \ \ \ \ \
-\,
\frac{R_{xyyy}}{R_x}
+
3\bigg(\!\frac{R_y}{R_x}\!\bigg)
\frac{R_{xxyy}}{R_x}
-
3\bigg(\!\frac{R_y}{R_x}\!\bigg)^2
\frac{R_{xxxy}}{R_x}
+
\bigg(\!\frac{R_y}{R_x}\!\bigg)^3
\frac{R_{xxxx}}{R_x}
\bigg],
\endaligned
\]
et bien sûr, la différentiation du membre de droite s'obtient
en échangeant $x \longleftrightarrow y$, avec
un signe <<\,$-$\,>> global:
\[
\aligned
&
-\,
\frac{x''''}{R_y}
-
\frac{x'''x'}{R_y}
\bigg[
-\,4\,\frac{R_{yx}}{R_y}
+
4\bigg(\!\frac{R_x}{R_y}\!\bigg)
\frac{R_{yy}}{R_y}
\bigg]
-
\frac{x''x''}{R_y}
\bigg[
-\,3\,\frac{R_{yx}}{R_y}
+
3\bigg(\!\frac{R_x}{R_y}\!\bigg)
\frac{R_{yy}}{R_y}
\bigg]
\,-
\\
&
-\,
\frac{x''(x')^2}{R_y}
\bigg[
3\,\frac{R_{yx}}{R_y}\frac{R_{yx}}{R_y}
+
3\,\frac{R_{xx}}{R_y}\frac{R_{yy}}{R_y}
-
30\bigg(\!\frac{R_x}{R_y}\!\bigg)
\frac{R_{yy}}{R_y}\frac{R_{yx}}{R_y}
+
15\,\bigg(\!\frac{R_x}{R_y}\!\bigg)^2
\frac{R_{yy}}{R_y}\frac{R_{yy}}{R_y}
\,-
\\
&
\ \ \ \ \ \ \ \ \ \ \ \ \ \ \ \
-\,
6\,\frac{R_{yxx}}{R_y}
+
12\bigg(\!\frac{R_x}{R_y}\!\bigg)
\frac{R_{yyx}}{R_y}
-
6\bigg(\!\frac{R_x}{R_y}\!\bigg)^2
\frac{R_{yyy}}{R_y}
\bigg]
\,-
\endaligned
\]
\[
\!\!\!\!\!\!\!\!\!\!\!\!\!\!\!\!\!\!\!\!\!\!\!\!\!\!
\aligned
&
-\,
\frac{(x')^4}{R_y}
\bigg[
-\,6\,\frac{R_{yy}}{R_y}\frac{R_{yx}}{R_y}\frac{R_{xx}}{R_y}
+
30\bigg(\!\frac{R_x}{R_y}\!\bigg)
\frac{R_{yy}}{R_y}\frac{R_{yx}}{R_y}\frac{R_{yx}}{R_y}
+
6\bigg(\!\frac{R_x}{R_y}\!\bigg)
\frac{R_{yy}}{R_y}\frac{R_{yy}}{R_y}\frac{R_{xx}}{R_y}
\,-
\\
&
\ \ \ \ \ \ \ \ \ \ \ \ \ \ \ \
-\,
45\bigg(\!\frac{R_x}{R_y}\!\bigg)^2
\frac{R_{yy}}{R_y}\frac{R_{yy}}{R_y}\frac{R_{yx}}{R_y}
+
15\bigg(\!\frac{R_x}{R_y}\!\bigg)^2
\frac{R_{yy}}{R_y}\frac{R_{yy}}{R_y}\frac{R_{yy}}{R_y}
\,+
\\
&
\ \ \ \ \ \ \ \ \ \ \ \ \ \ \ \
+
2\,\frac{R_{yx}}{R_y}\frac{R_{yxx}}{R_y}
+
2\,\frac{R_{xx}}{R_y}\frac{R_{yyx}}{R_y}
-
8\bigg(\!\frac{R_x}{R_y}\!\bigg)
\frac{R_{yy}}{R_y}\frac{R_{yxx}}{R_y}
-
14\bigg(\!\frac{R_x}{R_y}\!\bigg)
\frac{R_{yx}}{R_y}\frac{R_{yyx}}{R_y}
-
2\bigg(\!\frac{R_x}{R_y}\!\bigg)
\frac{R_{xx}}{R_y}
\frac{R_{yyy}}{R_y}
\,+
\\
&
\ \ \ \ \ \ \ \ \ \ \ \ \ \ \ \
+
18\bigg(\!\frac{R_x}{R_y}\!\bigg)^2\frac{R_{yy}}{R_y}
\frac{R_{yyx}}{R_y}
+
12\bigg(\!\frac{R_x}{R_y}\!\bigg)^2
\frac{R_{yx}}{R_y}\frac{R_{yyy}}{R_y}
-
10\bigg(\!\frac{R_x}{R_y}\!\bigg)^3
\frac{R_{yy}}{R_y}\frac{R_{yyy}}{R_y}
\,-
\\
&
\ \ \ \ \ \ \ \ \ \ \ \ \ \ \ \
-\,
\frac{R_{yxxx}}{R_y}
+
3\bigg(\!\frac{R_x}{R_y}\!\bigg)
\frac{R_{yyxx}}{R_y}
-
3\bigg(\!\frac{R_x}{R_y}\!\bigg)^2
\frac{R_{yyyx}}{R_y}
+
\bigg(\!\frac{R_x}{R_y}\!\bigg)^3
\frac{R_{yyyy}}{R_y}
\bigg].
\endaligned
\]

\section{\bf Annulations sur la droite à l'infini $\P_\infty^1$}
\label{annulations-infini}
\HEAD{\ref{annulations-infini}.~Annulations sur la droite à 
l'infini $\P_\infty^1$}{
Jo\"el Merker, Département de Mathématiques d'Orsay}

Comment ces différentielles de jets, holomorphes sur la partie affine
$X^1 \cap \P^2$ de la courbe, se comportent-elles aux points (isolés,
en nombre fini) $X^1 \cap \P_\infty^1$ de la courbe qui sont situés à
l'infini? Ces différentielles de jets y ont-elles des pôles?

\medskip

Dans la carte affine initiale $\{ T \neq 0\}$, {\sl l'équation affine
initiale} de $X^1 \cap {\sf U}_0$:
\[
0
=
R_0(x_0,y_0),
\] 
s'obtient en extrayant $T^d$ par
homogénéité:
\[
0
=
T^d\,
\underbrace{
{\rm R}\bigg(
1\colon\!
\frac{X}{T}
\colon\!
\frac{Y}{T}
\bigg)}_{
=:\,R_0(x_0,y_0)}.
\] 
De manière analogue, des divisions de l'équation homogène
par $X$ puis par $Y$:
\[
0
=
X^d\,
\underbrace{
{\rm R}\bigg(
\frac{T}{X}
\colon\!
1
\colon\!
\frac{Y}{X}
\bigg)}_{
=:\,R_1(x_1,y_1)}
\ \ \ \ \ \ \ \ \ \ \ \ \ \ \ \ \
\text{\rm et}
\ \ \ \ \ \ \ \ \ \ \ \ \ \ \ \ \
0
=
Y^d\,
\underbrace{
{\rm R}\bigg(
\frac{T}{Y}
\colon\!
\frac{X}{Y}
\colon\!
1
\bigg)}_{
=:\,R_2(x_2,y_2)},
\]
fournissent les deux autres équations affines.

\medskip

Sachant que les deux autres systèmes de coordonnées affines
sur $\P^2$: 
\[
(x_1,y_1)
\ \ \ \ \ \ \ \ \ \ \ \ \
\text{\rm et}
\ \ \ \ \ \ \ \ \ \ \ \ \
(x_2,y_2)
\]
sont liés aux coordonnées initiales $(x_0,y_0)$ comme suit:
\[
\aligned
(x_1,y_1)
=
\bigg(
\frac{1}{x_0},\,
\frac{y_0}{x_0}
\bigg)
&
\,\,\,\Longleftrightarrow\,\,\,
\bigg(
\frac{1}{x_1},\,
\frac{y_1}{x_1}
\bigg)
=
(x_0,y_0),
\\
(x_2,y_2)
=
\bigg(
\frac{x_0}{y_0},\,
\frac{1}{y_0}
\bigg)
&
\,\,\,\Longleftrightarrow\,\,\,
\bigg(
\frac{x_2}{y_2},\,
\frac{1}{y_2}
\bigg)
=
(x_0,y_0),
\endaligned
\]
l'homogénéité du polynôme ${\rm R} = {\rm R} 
(T \colon X \colon Y)$ montre par le calcul que:
\[
\aligned
R_1\big(x_1,y_1\big)
&
=
{\rm R}\bigg(
\frac{T}{X}
\colon\!
1
\colon\!
\frac{Y}{X}
\bigg)
=
\frac{T^d}{X^d}\,
{\rm R}\bigg(
1
\colon\!
\frac{X}{T}
\colon\!
\frac{Y}{T}
\bigg)
=
\frac{1}{(x_0)^d}\,
R\big(x_0,y_0\big),
\\
R_2\big(x_2,y_2\big)
&
=
{\rm R}\bigg(
\frac{T}{Y}
\colon\!
\frac{X}{Y}
\colon\!
1
\bigg)
=
\frac{T^d}{Y^d}\,
{\rm R}\bigg(
1
\colon\!
\frac{X}{T}
\colon\!
\frac{Y}{T}
\bigg)
=
\frac{1}{(y_0)^d}\,
R_0\big(x_0,y_0\big).
\endaligned
\]

Après normalisations géométriques élémentaires, dans
les coordonnées affines initiales abrégées:
\[
(x,y)
\equiv
(x_0,y_0),
\]
le polynôme définissant initial $R \equiv R_0$ s'écrit:
\[
R
=
x^d
+
y^d
+
R_d^*
+
R_{d-1}
+\cdots+
R_1
+
R_0,
\]
avec:
\[
R_d^*
=
\sum_{1\leqslant i\leqslant d-1}\,
\underbrace{r_{i,d-i}}_{\in\,\C}\,
x^i\,y^{d-i},
\]
et, pour $1 \leqslant \ell \leqslant d-1$, avec:
\[
\
R_\ell
=
\sum_{0\leqslant i\leqslant\ell}\,
\underbrace{r_{i,\ell-i}}_{\in\,\C}\,
x^i\,y^{\ell-i},
\]
le dernier terme étant constant:
\[
R_0
\in
\C.
\]
La présence de $x^d + y^d$\,\,---\,\,toujours 
arrangeable\,\,---\,\,assure que:
\[
\infty_x
\,\not\in\,
X^1\cap\P_\infty^1
\ \ \ \ \ \ \ \ \ \ \ \ \
\text{\rm et}
\ \ \ \ \ \ \ \ \ \ \ \ \
\infty_y
\,\not\in\,
X^1\cap\P_\infty^1.
\]
Par disposition géométrique déjà réalisée à l'avance, la courbe $X^1$
intersecte transversalement la droite à l'infini $\P_\infty^1$
en exactement $d$ points.

Il convient maintenant de ré-écrire les coordonnées affines initiales
avec l'indice qui notifie leur appartenance à l'ouvert ${\sf U}_0$:
\[
(x_0,y_0).
\]

\medskip\noindent{\bf Question de transfert.}
{\em Que devient alors, par exemple, la différentielle de jets:}
\[
\frac{y_0'}{R_{0,x_0}(x_0,y_0)}
\bigg\vert_{X^1}
\]
{\em à travers le changement de carte affine:}
\[
(x_0,y_0)
=
\bigg(
\frac{x_2}{y_2},\,\frac{1}{y_2}
\bigg)\,?
\]

\medskip
Dans l'ouvert ${\sf U}_2$, la nouvelle équation affine de la courbe:
\[
0
=
R_2(x_2,y_2)
\]
est donnée en termes du nouveau polynôme:
\[
\aligned
R_2(x_2,y_2)
:=
&\,
(y_2)^d\,
R_0
\bigg(
\frac{x_2}{y_2},\,\frac{1}{y_2}
\bigg)
\\
=
&\,
\big(y_2\big)^d\,
\bigg[
\bigg(\!\frac{x_2}{y_2}\!\bigg)^d
+
\bigg(\!\frac{1}{y_2}\!\bigg)^d
+
\sum_{1\leqslant i\leqslant d-1}\,
r_{i,d-i}\,
\bigg(\!\frac{x_2}{y_2}\!\bigg)^i
\bigg(\!\frac{1}{y_2}\!\bigg)^{d-i}
+
\\
&
\ \ \ \ \ \ \ \ \ \ \ \ \ \ \ \ \ \ \ \ \ \ \ 
+
R_{d-1}
\bigg(
\frac{x_2}{y_2},\,\frac{1}{y_2}
\bigg)
+\cdots+
R_1
\bigg(
\frac{x_2}{y_2},\,\frac{1}{y_2}
\bigg)
+
R_0
\bigg].
\endaligned
\]

Or, chaque $R_\ell$ est homogène de degré $\ell$:
\[
\big(y_2\big)^d\,
R_\ell
\bigg(
\frac{x_2}{y_2},\,\frac{1}{y_2}
\bigg)
\,=\,
{\rm O}
\Big(
\big(y_2\big)^{d-\ell}
\Big)
\ \ \ \ \ \ \ \ \ \ \ \ \
{\scriptstyle{(0\,\leqslant\,\ell\,\leqslant\,d\,-\,1)}},
\]
donc:
\[
R_2(x_2,y_2)
=
(x_2)^d
+
1
+
\sum_{1\leqslant i\leqslant d-1}\,
r_{i,d-i}\,(x_2)^i
+
{\rm O}(y_2).
\]

Dans ces nouvelles coordonnées affines $(x_2,y_2)$:
\[
\P_\infty^1
\big\backslash\{\infty_x\}
\,=\,
\big\{y_2=0\big\},
\]
et aussi le point:
\[
\infty_y
=
\{x_2=0,\,y_2=0\}
\]
est devenu l'origine sur cette droite. En restriction à cette droite
affine $\{ y_2 = 0\}$:
\[
R_2(x_2,0)
=
(x_2)^d
+
1
+
\sum_{1\leqslant i\leqslant d-1}\,r_{i,d-i}\,(x_2)^i,
\]
ce qui montre bien que $\infty_y \not\in X^1$.

\medskip

De plus, les $d$ racines de $R_2(x_2, 0) = 0$: 
\[
\big(\underline{x}_2^1,0\big),
\,\dots\dots,\,
\big(\underline{x}_2^d,0\big),
\]
ensemble des points $X^1 \cap \P_\infty^1$, 
sont mutuellement distinctes, toutes de multiplicité
$1$, ce qui implique (transversalité, non-tangentialité):
\[
0
\neq
\frac{\partial R_{2,x_2}}{\partial x_2}
\big(\underline{x}_2^1,0\big),
\,\,\dots\dots\dots,\,\,
0
\neq
\frac{\partial R_{2,x_2}}{\partial x_2}
\big(\underline{x}_2^d,0\big).
\] 

Ensuite, à travers le $(1/y)$-changement de carte, une différentiation de:
\[
y_0
=
\frac{1}{y_2},
\]
donne:
\[
\boxed{\,
y_0'
=
-\,\frac{y_2'}{(y_2)^2},\,}
\]
ce qui traite le transfert du numérateur de la différentielle
de jets $\frac{y_0'}{R_{0,x_0}(x_0,y_0)}$, tandis qu'une
différentiation de l'identité équivalente
à celle définissant $R_2$ vue plus haut:
\[
R_0(x_0,y_0)
\,\equiv\,
(y_0)^d\,
R_2
\bigg(
\frac{x_0}{y_0},\,\frac{1}{y_0}
\bigg)
\]
par rapport à $x_0$ donne:
\[
R_{0,x_0}(x_0,y_0)
\,=\,
\big(y_0\big)^{d-1}\,
R_{2,x_2}
\bigg(
\frac{x_0}{y_0},\,\frac{1}{y_0}
\bigg),
\]
à savoir:
\[
\boxed{\,
R_{0,x_0}(x_0,y_0)
\,=\,
\frac{1}{(y_2)^{d-1}}\,
R_{2,x_2}(x_2,y_2).\,}
\]

Tout cela mis ensemble montre comment se transfère la différentielle de jets
en question:
\[
\aligned
\frac{y_0'}{R_{0,x_0}(x_0,y_0)}
&
=
\frac{-\,\frac{y_2'}{(y_2)^2}}{
\frac{1}{(y_2)^{d-1}}\,R_{2,x_2}(x_2,y_2)}
\\
&
=
-\,y_2'\,
\frac{(y_2)^{d-3}}{R_{2,x_2}(x_2,y_2)}.
\endaligned
\]

Puisque la courbe $X^1$ intersecte transversalement 
la droite $\{ y_2 = 0\}$\,\,---\,\,
en laquelle s'est transformé $\P_\infty^1 \backslash 
\{ \infty_x\}$\,\,---\,\,exactement aux
$d$ points mutuellement distincts:
\[
\big(\underline{x}_2^1,0\big),
\,\dots\dots,\,
\big(\underline{x}_2^d,0\big),
\]
les valeurs en ces points du dénominateur:
\[
\frac{1}{
R_{2,x_2}(x_2,y_2)}
\]
sont non nulles, donc localement holomorphes en tant que
fonctions de $(x_2, y_2)$, et alors, pourvu seulement que:
\[
\boxed{\,d\,\geqslant\,3,\,}
\]
le tout est vraiment holomorphe. Plus encore, et sous
forme de résumé synthétique auto-contenu:

\begin{Theoreme}
{\small\sf (XIX\textsuperscript{ème} siècle, {\em cf.}
\cite{Griffiths-1989})}
Sur une courbe algébrique projective lisse quelconque:
\[
X^1
\,\subset\,
\P^2(\C)
\]
de degré:
\[
d\geqslant 4,
\]
il existe toujours un système de coordonnées affines:
\[
(x,y)
\,\in\,
\C^2
\subset
\P^2
\]
dans lesquelles la disposition géométrique de la courbe:
\[
X^1
=
\underbrace{\big(X^1\cap\C^2\big)}_{
\text{\sf partie affine}}
\,\bigcup\,
\underbrace{
\big(X^1\cap\P_\infty^1\big)}_{
\text{\sf points à l'infini}}
\]
est telle que l'intersection:
\[
X^1
\cap
\P_\infty^1
=:
X_\infty^0
\]
consiste en exactement $d$ points distincts deux à deux et distincts
des deux points $\infty_x$ et $\infty_y$ situés à l'infini sur l'axe des $x$
et à l'infini 
sur l'axe des $y$, en lesquels la droite projective $\P_\infty^1$
est non tangente à la courbe $X^1$, et dans de telles circonstantes,
si l'équation affine s'exprime comme lieu des zéros:
\[
\Big\{
(x,y)\in\C^2\colon\,
R(x,y)=0
\Big\}
\]
d'un certain polynôme $R = R(x, y)$ de degré $d \geqslant 4$, 
l'expression:
\[
\boxed{\,
{\sf J}_R^1
\,:=\,
\left\{
\aligned
&
\ \ \ \ \ \,
\frac{y'}{R_x(x,y)}
\bigg\vert_{X^1\,\cap\,\C^2\,\cap\,\{R_x\neq 0\}},
\\
&
-\,
\frac{x'}{R_y(x,y)}
\bigg\vert_{X^1\,\cap\,\C^2\,\cap\,\{R_y\neq 0\}},\,
\\
&
\ \ \ \ \ \ \ \ \ \ \ 
0
\ \ \ \ \ \ \ \ \ 
\text{\rm sur}\ \
X^1\cap\P_\infty^1,
\endaligned
\right.}
\]
définit une \underline{\sl différentielle génératrice} de $1$-jets
{\em holomorphe} sur $X^1$ tout entier qui s'annule identiquement sur
le diviseur ample $X^1 \cap \P_\infty^1$.
\end{Theoreme}

\proof
Il suffit d'observer d'abord que le transfert à l'infini de $\frac{
y'}{ R_x}$, à savoir:
\[
-\,y_2'\,
\frac{(y_2)^{d-3}}{
R_{2,x_2}(x_2,y_2)}
\,=\,
y_2'\,{\rm O}(y_2)
\]
s'annule identiquement sur $\{ y_2 = 0\}$. Ensuite, il reste à
mentionner que le fibré holomorphe des $1$-jets de courbes holomorphes
$\D \longrightarrow X^1$:
\[
\xymatrix{
J^1(\D,X^1) \dto^\pi
\\
X^1,
}
\]
peut être saisi, au-dessus de $X^1 \cap \C^2$, dans deux collections
de systèmes de cartes holomorphes locales contenues
respectivement dans:
\[
X^1\cap\C^2\cap\{R_x\neq 0\},
\]
et dans:
\[
X^1\cap\C^2\cap\{R_y\neq 0\},
\]
par simples projections sur l'axe des $y$ et sur l'axe
des $y$, respectivement.
\endproof

La manière dont une telle {\em différentielle génératrice}
engendre toutes les sections holomorphes globales
de ce fibré $J^1 ( \D, X^1)$ sera explicitée ultérieurement.

\medskip

Mais auparavant, il s'agit d'explorer comment se transfèrent à l'infini
les autres différentielles génératrices construites précédemment.

\medskip

Mieux vaut étudier tout d'abord un exemple. Pour les jets d'ordre
$\kappa = 3$, le membre de gauche de l'équation génératrice:
\[
\aligned
&
\frac{y'''}{R_x}
+
\underline{\frac{y''y'}{R_x}
\bigg[
-\,3\,
\frac{R_{xy}}{R_x}}
+
3\bigg(\!\frac{R_y}{R_x}\!\bigg)
\frac{R_{xx}}{R_x}
\bigg]
\,+
\\
&
\ \ \ \ \ \
+
\frac{(y')^3}{R_x}
\bigg[
-\,
6\,\bigg(\!\frac{R_y}{R_x}\!\bigg)
\frac{R_{xx}}{R_x}
\frac{R_{xy}}{R_x}
+
3\bigg(\!\frac{R_y}{R_x}\!\bigg)^2
\frac{R_{xx}}{R_x}
\frac{R_{xx}}{R_x}
\,-
\\
&
\ \ \ \ \ \ \ \ \ \ \ \ \ \ \ \ \ \ \ \ \ 
-\,
\frac{R_{xyy}}{R_x}
+
2\,\bigg(\!\frac{R_y}{R_x}\!\bigg)
\frac{R_{xxy}}{R_x}
-
\bigg(\!\frac{R_y}{R_x}\!\bigg)^2
\frac{R_{xxx}}{R_x}
\bigg]
\,=\,
\\
&
\,=\,
-\,
\frac{x'''}{R_y}
-
\frac{x''x'}{R_y}
\bigg[
-\,3\,
\frac{R_{xy}}{R_y}
+
3\bigg(\!\frac{R_x}{R_y}\!\bigg)
\frac{R_{yy}}{R_y}
\bigg]
\,+
\\
&
\ \ \ \ \ \
+
\frac{(x')^3}{R_y}
\bigg[
-\,
6\,\bigg(\!\frac{R_x}{R_y}\!\bigg)
\frac{R_{yy}}{R_y}
\frac{R_{xy}}{R_y}
+
3\bigg(\!\frac{R_x}{R_y}\!\bigg)^2
\frac{R_{yy}}{R_y}
\frac{R_{yy}}{R_y}
\,-
\\
&
\ \ \ \ \ \ \ \ \ \ \ \ \ \ \ \ \ \ \ \ \ 
-\,
\frac{R_{xxy}}{R_y}
+
2\,\bigg(\!\frac{R_x}{R_y}\!\bigg)
\frac{R_{xyy}}{R_y}
-
\bigg(\!\frac{R_x}{R_y}\!\bigg)^2
\frac{R_{yyy}}{R_y}
\bigg],
\endaligned
\]
comporte par exemple le terme (souligné, en négligeant
le coefficient $-3$):
\[
\frac{y''y'}{R_x}\,
\frac{R_{xy}}{R_x}.
\]

\medskip\noindent{\bf Question.}
{\em Comment ce terme
$\frac{ y'' y'}{ R_x}\, \frac{ R_{ xy}}{R_x}$ se transfère-t-il à
l'infini?}

\medskip
En tout cas, ce qui vient d'être vu a montré que:
\[
d-3
=
\ordreinfini
\bigg(
\frac{y'}{R_x(x,y)}
\bigg).
\]
Qu'en est-il, alors, de: 
\[
\frac{y''\,y'}{ R_x}\,\frac{R_{xy}}{R_x}
=
\frac{y_0''\,y_0'}{R_{0,x_0}(x_0,y_0)}\,
\frac{R_{0,x_0y_0}}{R_{0,x_0}(x_0,y_0)}\,?
\]

C'est assez simple. En revenant donc, comme il se doit, à la notation
$(x_0, y_0)$, une différentiation supplémentaire de:
\[
y_0'
=
-\,\frac{y_2'}{(y_2)^2}
\]
donne instantanément:
\[
y_0''
=
-\,\frac{y_2''}{(y_2)^2}
+
2\,\frac{y_2'\,y_2'}{(y_2)^3}.
\]
Ce qui vient d'être vu a aussi montré que:
\[
\frac{1}{R_{0,x_0}(x_0,y_0)}
\,=\,
\frac{(y_2)^{d-1}}{R_{2,x_2}(x_2,y_2)}.
\]
Par ailleurs:
\[
R_{0,x_0y_0}\big(x_0,y_0\big)
\,=\,
R_{0,x_0y_0}
\bigg(
\frac{x_2}{y_2},\frac{1}{y_2}
\bigg).
\]
Donc:
\[
\frac{y_0''\,y_0'}{R_{0,x_0}(x_0,y_0)}\,
\frac{R_{0,x_0y_0}(x_0,y_0)}{R_{0,x_0}(x_0,y_0)}
\,=\,
\frac{
\left(
-\,\frac{y_2''}{(y_2)^2}
+
2\,
\frac{y_2'\,y_2'}{(y_2)^3}
\right)\,
\left(
-\,\frac{y_2'}{(y_2)^2}
\right)}{
\frac{1}{(y_2)^{d-1}}\,
R_{2,x_2}(x_2,y_2)}
\,\,
\frac{
R_{0,x_0y_0}
\left(
\frac{x_2}{y_2},\frac{1}{y_2}
\right)
}{
\frac{1}{(y_2)^{d-1}}\,
R_{2,x_2}(x_2,y_2)}.
\]

Maintenant, puisque le polynôme $R_0$ est de degré exactement
égal à $d$, sa dérivée d'ordre deux $R_{0, x_0y_0}$ est de
degré exactement égal à $d-2$.

Donc après réduction au même dénominateur:
\[
\aligned
R_{0,x_0y_0}
\bigg(
\frac{x_2}{y_2},\frac{1}{y_2}
\bigg)
&
\,=\,
\frac{1}{(y_2)^{d-2}}\,
\mathmotsf{polynôme}(x_2,y_2)
\\
&
\,
=:
\frac{1}{(y_2)^{d-2}}\,
S_{1,1}(x_2,y_2).
\endaligned
\]

Au total:
\[
\aligned
\frac{y_0''\,y_0'}{R_{0,x_0}(x_0,y_0)}\,
\frac{R_{0,x_0y_0}(x_0,y_0)}{R_{0,x_0}(x_0,y_0)}
&
\,=\,
\frac{(y_2)^{d-1-2-3}
\big(
y_2''y_2'y_2
-
2\,(y_2')^3
\big)}{
R_{2,x_2}(x_2,y_2)}
\cdot
\frac{\boxed{(y_2)^1}\,S_{1,1}(x_2,y_2)}{
R_{2,x_2}(x_2,y_2)}
\\
&
\,=\,
\frac{(y_2)^{d-5}\,\big(y_2''y_2'y_2-2\,(y_2')^3\big)
\cdot
S_{1,1}(x_2,y_2)}{
R_{2,x_2}(x_2,y_2)}
\endaligned
\]
Ici, le terme encadré:
\[
\boxed{
(y_2)^1}
=
\frac{
\frac{1}{(y_2)^{d-2}}
}{
\frac{1}{(y_2)^{d-1}}
}
\]
fait apparaître un exposant {\em positif}
qu'il convient de conceptualiser comme:
\[
\aligned
\infinibonus
\bigg(
\frac{R_{0,x_0y_0}}{R_{0,x_0}}
\bigg)
:=
&\,
\ordreinfini
\bigg(
\frac{R_{0,x_0y_0}}{R_{0,x_0}}
\bigg)
\\
=
&\,
1,
\endaligned
\]
ce qui conclut que:
\[
\aligned
\ordreinfini
\bigg(
\frac{y''y'}{R_x}\,
\frac{R_{xy}}{R_x}
\bigg)
&
=
d-6+1
&
=
d-5.
\endaligned
\]

\begin{Lemme}
Tous les termes dans la première ligne de l'équation génératrice
ci-dessus au niveau des jets d'ordre $3$, par exemple les cinq
termes:
\[
\frac{y'''}{R_x},
\ \ \ \ \ \ \ \ \ 
\underbrace{\frac{y''y'}{R_x}\,
\frac{R_{xy}}{R_x}}_{
\text{\sf vient d'être vu}},
\ \ \ \ \ \ \ \ \ 
\frac{y''y'}{R_x}
\bigg(\!\frac{R_y}{R_x}\!\bigg)
\frac{R_{xx}}{R_x},
\ \ \ \ \ \ \ \ \ 
\frac{(y')^3}{R_x}
\bigg(\!\frac{R_y}{R_x}\!\bigg)
\frac{R_{xx}}{R_x}
\frac{R_{xy}}{R_x},
\ \ \ \ \ \ \ \ \ 
\frac{(y')^3}{R_x}
\frac{R_{xyy}}{R_x},
\]
ont un ordre à l'infini égal à:
\[
d-5.
\]
\end{Lemme}

\proof
En effet, trois différentiations de $y_0 = \frac{ 1}{ y_2}$ donnent:
\[
\aligned
y_0'
&
=
-\,\frac{y_2'}{(y_2)^2},
\\
y_0''
&
=
-\,\frac{y_2''}{(y_2)^2}
+
2\,
\frac{y_2'y_2'}{(y_2)^3},
\\
y_0'''
&
=
-\,\frac{y_2'''}{(y_2)^2}
+
4\,
\frac{y_2''y_2'}{(y_2)^3}
-
6\,\frac{(y_2')^3}{(y_2)^4},
\endaligned
\]
d'où:
\[
\aligned
\ordreinfini
\big(y_0'\big)
&
=
-\,2,
\\
\ordreinfini
\big(y_0''\big)
&
=
-\,3,
\\
\ordreinfini
\big(y_0'''\big)
&
=
-\,4.
\endaligned
\]

Par ailleurs, sachant que:
\[
\ordreinfini
\bigg(
\frac{1}{R_x}
\bigg)
=
d-1,
\]
il vient:
\[
\aligned
\ordreinfini
\bigg(
\frac{y'''}{R_x}
\bigg)
&
=
-\,4+d-1
\\
&
=
d-5.
\endaligned
\]

De même, sachant que:
\[
\ordreinfini
\big(R_y\big)
=
\ordreinfini
\big(R_x\big)
=
d-1,
\]
il vient:
\[
\ordreinfini
\bigg(
\frac{R_y}{R_x}
\bigg)
=
0,
\]
puis sachant que:
\[
\aligned
\ordreinfini
\bigg(
\frac{R_{xx}}{R_x}
\bigg)
=
\underbrace{-\,
(d-2)+(d-1)}_{
\infinibonus\,=\,1},
\endaligned
\]
il vient:
\[
\aligned
\ordreinfini
\bigg(
\frac{y''y'}{R_x}
\bigg(\!\frac{R_y}{R_x}\!\bigg)
\frac{R_{xx}}{R_x}
\bigg)
&
=
-\,3-2+(d-1)+0+1
\\
&
=
d-5.
\endaligned
\]

Le quatrième monôme sélectionné par l'énoncé se traite
de la même manière.

Quant au dernier monôme, il a le même:
\[
\aligned
\ordreinfini
\bigg(
\frac{(y')^3}{R_x}\,
\frac{R_{xyy}}{R_x}
\bigg)
&
=
-\,4+(d-1)
+
2
\\
&
=
d-5,
\endaligned
\]
puisque (exercice):
\[
\aligned
\infinibonus
\bigg(
\frac{R_{xyy}}{R_x}
\bigg)
&
=
\ordreinfini
\big(
R_{xyy}
\big)
-
\ordreinfini
\big(
R_{x}
\big)
\\
&
=
-\,(d-3)
+
(d-1)
\\
&
=
2,
\endaligned
\]
ce qui conclut pour les cinq exemples de monômes, et tous les autres
monômes s'avèrent posséder le {\em même} ordre $d-5$ à l'infini
(exercice visuel).
\endproof

\begin{Lemme}
\'Etant donné un ordre de jets quelconque $\kappa \geqslant 1$, 
pour tout ordre intermédiaire $1 \leqslant \lambda \leqslant \kappa$,
la dérivée d'ordre $\lambda$ de:
\[
y_0
=
\frac{1}{y_2}
\]
produit une formule du type:
\[
\aligned
y_0^{(\lambda)}
&
=
-\,\frac{y_2^{(\lambda)}}{(y_2)^2}
+\cdots+
(-1)^\lambda\,\lambda!\,
\frac{(y_2')^\lambda}{(y_2)^{\lambda+1}}
\\
&
=
\frac{\text{\scriptsize\sf polynôme}_\Z^\lambda
\big(y_2,y_2',\dots,y_2^{(\lambda)}\big)}{
(y_2)^{\lambda+1}},
\endaligned
\]
ce dernier polynôme à coefficients dans
$\Z$ ayant la propriété que tous ses monômes:
\[
\big(y_2\big)^{\nu_0}
\big(y_2'\big)^{\nu_1}
\cdots
\big(y_2^{(\lambda)}\big)^{\nu_\lambda}
\]
ont uniformément le même nombre total de `primes':
\[
1\cdot\nu_1
+\cdots+
\lambda\cdot\nu_\lambda
=
\lambda.
\]
\end{Lemme}

\proof
Laissée au lecteur, il s'agit d'une récurrence peu exigeante dans
laquelle il est décidé de ne pas expliciter ces polynômes.
\endproof

Ainsi généralement:
\[
\boxed{\,
\ordreinfini
\big(
y_0^{(\lambda)}
\big)
=
-\,\lambda-1.\,}
\]

Plus encore, pour tous entiers-exposants:
\[
\mu_1\geqslant 0,
\ \ \
\mu_2\geqslant 0,
\ \ \
\dots\dots,
\ \ \
\mu_{\kappa-1}\geqslant 0,
\]
puisque:
\[
\!\!\!\!\!\!\!\!\!\!\!\!\!\!\!\!\!\!\!\!
\aligned
&
\big(y_0'\big)^{\mu_1}
\big(y_0''\big)^{\mu_2}
\cdots
\big(y_0^{(\kappa-1)}\big)^{\mu_{\kappa-1}}
=
\\
&
\ \ \ \ \
=
\bigg(\!
-\,\frac{y_2'}{(y_2)^2}
\!\bigg)^{\mu_1}
\bigg(\!
-\,\frac{y_2''}{(y_2)^2}
+
2\,\frac{y_2'y_2'}{(y_2)^3}
\!\bigg)^{\mu_2}
\cdots\cdots
\bigg(\!
-\,\frac{y_2^{(\kappa-1)}}{(y_2)^2}
+\cdots+
(-1)^{\kappa-1}
(\kappa-1)!
\frac{(y_2')^{\kappa-1}}{(y_2)^\kappa}
\!\bigg)^{\mu_{\kappa-1}}
\\
&
\ \ \ \ \
=
\frac{\text{\scriptsize\sf Polynôme}
\big(y_2,y_2',\dots,y_2^{(\kappa-1)}\big)}{
(y_2)^{2\mu_1+3\mu_2+\cdots+\kappa\mu_{\kappa-1}}},
\endaligned
\]
il vient:
\[
\boxed{\,
\ordreinfini
\Big(
\big(y_0'\big)^{\mu_1}
\big(y_0''\big)^{\mu_2}
\cdots\cdots
\big(y_0^{(\kappa-1)}\big)^{\mu_{\kappa-1}}
\Big)
=
-\,
2\,\mu_1
-
3\,\mu_2
-\cdots-
\kappa\,\mu_{\kappa-1}.\,}
\]

\medskip

Maintenant, dans l'équation génératrice d'une différentielle
de jets d'ordre quelconque $\kappa \geqslant 1$ fixé:
\[
\!\!\!\!\!\!\!\!\!\!\!\!\!\!\!\!\!\!\!\!
\aligned
&
\frac{y^{(\kappa)}}{R_x}
+
\sum_{\mu_1+\cdots+(\kappa-1)\mu_{\kappa-1}=\kappa
\atop
\mu_1\geqslant 0,\dots,\mu_{\kappa-1}\geqslant 0}
\!\!\!\!\!
\frac{\big(y'\big)^{\mu_1}\cdots\big(y^{(\kappa-1)}\big)^{\mu_{\kappa-1}}}{
R_x}\,
\mathcal{J}_{\mu_1,\dots,\mu_{\kappa-1}}^\kappa
\left(
\frac{R_y}{R_x},\,\,
\bigg(
\frac{R_{x^i y^j}}{R_x}
\bigg)_{2\leqslant i+j\leqslant
\atop
\leqslant
-1+\mu_1+\cdots+\mu_{\kappa-1}}
\right)
\,=\,
\\
&
\ \ \ \ \
\,=\,
-\,
\frac{x^{(\kappa)}}{R_y}
-
\sum_{\mu_1+\cdots+(\kappa-1)\mu_{\kappa-1}=\kappa
\atop
\mu_1\geqslant 0,\dots,\mu_{\kappa-1}\geqslant 0}
\!\!\!\!\!
\frac{\big(x'\big)^{\mu_1}\cdots\big(x^{(\kappa-1)}\big)^{\mu_{\kappa-1}}}{
R_y}\,
\mathcal{J}_{\mu_1,\dots,\mu_{\kappa-1}}^\kappa
\left(
\frac{R_x}{R_y},\,\,
\bigg(
\frac{R_{y^ix^j}}{R_y}
\bigg)_{2\leqslant i+j\leqslant
\atop
\leqslant
-1+\mu_1+\cdots+\mu_{\kappa-1}}
\right),
\endaligned
\]
une information supplémentaire est aisément vérifiable par
récurrence.

\begin{Lemme}
Les polynômes apparaissant:
\[
\mathcal{J}_{\mu_1,\dots,\mu_{\kappa-1}}^\kappa
=
\mathcal{J}_{\mu_1,\dots,\mu_{\kappa-1}}^\kappa
\bigg(
{\sf R}_{0,1},\,\,
\Big(
{\sf R}_{i,j}
\Big)_{2\leqslant i+j\leqslant
-1+\mu_1+\cdots+\mu_{\kappa-1}}
\bigg),
\]
à coefficients dans $\Z$, ont tous leurs monômes de la forme:
\[
\big({\sf R}_{0,1}\big)^h\,
{\sf R}_{i_1,j_1}\,
\cdots
{\sf R}_{i_\nu,j_\nu}
\ \ \ \ \ \ \ \ \ \ \ \ \
{\scriptstyle{
(h\,\geqslant\,0;\,\,
2\,\leqslant\,
i_1+j_1,\,\dots,\,i_\nu+j_\nu
\,\leqslant\,
-\,1+\mu_1+\cdots+\mu_{\kappa-1})}},
\]
avec constamment:
\[
\big(i_1+j_1-1\big)
+\cdots+
\big(i_\nu+j_\nu-1\big)
=
-\,1
+
\mu_1
+\cdots+
\mu_{\kappa-1},
\]
de telle sorte que:
\[
\infinibonus
\bigg(
\bigg(\!\frac{R_y}{R_x}\!\bigg)^{\!h}\,
\frac{R_{x^{i_1}y^{i_1}}}{R_x}
\cdots
\frac{R_{x^{i_\nu}y^{i_\nu}}}{R_x}
\bigg)
=
-\,1
+
\mu_1
+\cdots+
\mu_{\kappa-1}.
\]
\end{Lemme}

\proof
Un examen des expressions complètes à l'ordre $\kappa = 2, 3, 4$
explicitées plus haut confirme cette constance combinatoire, dont la
démonstration, aisée, est laissée au lecteur.

Ensuite, à travers le changement de carte affine:
\[
(x_0,y_0)
\,\longmapsto\,
\bigg(
\frac{x_0}{y_0},
\frac{1}{y_0}
\bigg)
=
(x_2,y_2),
\]
chaque tel monôme de fractions de dérivées partielles\,\,---\,\,en
négligeant $\big( \frac{ R_y}{R_x} \big)^h$
qui ne change rien\,\,---\,\,se transforme comme:
\[
\!\!\!\!\!\!\!\!\!\!\!\!\!\!\!\!\!\!\!\!
\!\!\!\!\!\!\!\!\!\!\!\!\!\!\!\!\!\!\!\!
\aligned
\frac{R_{0,x_0^{i_1}y_0^{j_1}}(x_0,y_0)}{
R_{0,x_0}(x_0,y_0)}
\,\cdots\cdots\,
\frac{R_{0,x_0^{i_\nu}y_0^{j_\nu}}(x_0,y_0)}{
R_{0,x_0}(x_0,y_0)}
&
=
\frac{R_{0,x_0^{i_1}y_0^{j_1}}
\big(\frac{x_2}{y_2},\frac{1}{y_2}\big)
}{
\frac{1}{(y_2)^{d-1}}\,
R_{2,x_2}(x_2,y_2)}
\,\cdots\cdots\,
\frac{R_{0,x_0^{i_\nu}y_0^{j_\nu}}
\big(\frac{x_2}{y_2},\frac{1}{y_2}\big)
}{
\frac{1}{(y_2)^{d-1}}\,
R_{2,x_2}(x_2,y_2)}
\\
&
=
\frac{
\frac{1}{(y_2)^{d-i_1-j_1}}\,
S_{i_1,j_1}(x_2,y_2)
}{
\frac{1}{(y_2)^{d-1}}\,
R_{2,x_2}(x_2,y_2)}
\,\cdots\cdots\,
\frac{
\frac{1}{(y_2)^{d-i_\nu-j_\nu}}\,
S_{i_\nu,j_\nu}(x_2,y_2)
}{
\frac{1}{(y_2)^{d-1}}\,
R_{2,x_2}(x_2,y_2)}
\\
&
=
\frac{(y_2)^{i_1+j_1-1}\,S_{i_1,j_1}(x_2,y_2)}{
R_{2,x_2}(x_2,y_2)}
\,\cdots\cdots\,
\frac{(y_2)^{i_\nu+j_\nu-1}\,S_{i_\nu,j_\nu}(x_2,y_2)}{
R_{2,x_2}(x_2,y_2)}
\\
&
=
\frac{(y_2)^{i_1+j_1-1+\cdots+i_\nu+j_\nu-1}}{
\big[R_{2,x_2}(x_2,y_2)\big]^\nu}\,
\text{\footnotesize\sf polynôme}
(x_2,y_2),
\endaligned
\]
ce qui conclut.
\endproof

\begin{Proposition}
L'ordre à l'infini de chaque monôme de chacun des deux
membres de l'équation génératrice d'une différentielle
de jets d'ordre quelconque $\kappa \geqslant 1$ holomorphe
sur la partie affine $\C^2 \cap X^1$ de la courbe:
\[
\!\!\!\!\!\!\!\!\!\!\!\!\!\!\!\!\!\!\!\!
\!\!\!\!\!\!\!\!\!\!\!
\aligned
&
\frac{y^{(\kappa)}}{R_x}
+
\sum_{\mu_1+\cdots+(\kappa-1)\mu_{\kappa-1}=\kappa
\atop
\mu_1\geqslant 0,\dots,\mu_{\kappa-1}\geqslant 0}
\!\!\!\!\!
\frac{\big(y'\big)^{\mu_1}\cdots\big(y^{(\kappa-1)}\big)^{\mu_{\kappa-1}}}{
R_x}\,
\mathcal{J}_{\mu_1,\dots,\mu_{\kappa-1}}^\kappa
\left(
\frac{R_y}{R_x},\,\,
\bigg(
\frac{R_{x^i y^j}}{R_x}
\bigg)_{2\leqslant i+j\leqslant
\atop
\leqslant
-1+\mu_1+\cdots+\mu_{\kappa-1}}
\right)
\,=\,
\\
&
\ \ \ \ \
\,=\,
-\,
\frac{x^{(\kappa)}}{R_y}
-
\sum_{\mu_1+\cdots+(\kappa-1)\mu_{\kappa-1}=\kappa
\atop
\mu_1\geqslant 0,\dots,\mu_{\kappa-1}\geqslant 0}
\!\!\!\!\!
\frac{\big(x'\big)^{\mu_1}\cdots\big(x^{(\kappa-1)}\big)^{\mu_{\kappa-1}}}{
R_y}\,
\mathcal{J}_{\mu_1,\dots,\mu_{\kappa-1}}^\kappa
\left(
\frac{R_x}{R_y},\,\,
\bigg(
\frac{R_{y^ix^j}}{R_y}
\bigg)_{2\leqslant i+j\leqslant
\atop
\leqslant
-1+\mu_1+\cdots+\mu_{\kappa-1}}
\right),
\endaligned
\]
vaut constamment:
\[
\boxed{\,
d-\kappa-2.\,}
\]
\end{Proposition}

\proof
En effet, chaque monôme du membre de gauche de cette équation
parfaitement symétrique s'écrit:
\[
\boxed{\,
\frac{
\big(y_0'\big)^{\mu_1}
\big(y_0''\big)^{\mu_2}
\cdots\cdots
\big(y_0^{(\kappa-1)}\big)^{\mu_{\kappa-1}}}{
R_x}
\cdot
\bigg(\!\frac{R_y}{R_x}\!\bigg)^{\!h}\,
\frac{R_{x^{i_1}y^{i_1}}}{R_x}
\cdots
\frac{R_{x^{i_\nu}y^{i_\nu}}}{R_x},\,}
\]
et il vient d'être constaté à l'instant que:
\[
\ordreinfini
\bigg(
\frac{
\big(y_0'\big)^{\mu_1}
\big(y_0''\big)^{\mu_2}
\cdots\cdots
\big(y_0^{(\kappa-1)}\big)^{\mu_{\kappa-1}}}{
R_x}
\bigg)
=
-\,2\mu_1-3\mu_2
-\cdots-
\kappa\mu_{\kappa-1}
+
(d-1)
\]
et que:
\[
\aligned
\ordreinfini
\bigg(
\bigg(\!\frac{R_y}{R_x}\!\bigg)^{\!h}\,
\frac{R_{x^{i_1}y^{j_1}}}{R_x}
\cdots
\frac{R_{x^{i_\nu}y^{j_\nu}}}{R_x}
\bigg)
&
=
i_1+j_1-1
+\cdots+
i_\nu+j_\nu-1
\\
&
=
-\,1
+
\mu_1
+\cdots+
\mu_{\kappa-1},
\endaligned
\]
donc une simple addition donne la valeur:
\[
-\,2\mu_1-3\mu_2
-\cdots-
\kappa\mu_{\kappa-1}
+
(d-1)
-
1
+
\mu_1
+\cdots+
\mu_{\kappa-1}
\]
qui, en tenant compte de la contrainte de sommation:
\[
\sum_{\mu_1+\cdots+(\kappa-1)\mu_{\kappa-1}=\kappa
\atop
\mu_1\geqslant 0,\dots,\mu_{\kappa-1}\geqslant 0}
\]
vaut bien $d - \kappa - 2$.
\endproof

\begin{Theoreme}
\'Etant donné un ordre de jets arbitraire:
\[
\kappa
\geqslant
1,
\]
sur une courbe algébrique projective lisse quelconque:
\[
X^1
\,\subset\,
\P^2(\C)
\]
de degré:
\[
d
\,\geqslant\,
\kappa+3,
\]
saisie comme précédemment dans un système de coordonnées affines
adaptées:
\[
(x,y)
\,\in\,
\C^2
\subset
\P^2
\]
comme lieu des zéros:
\[
\Big\{
(x,y)\in\C^2\colon\,
R(x,y)=0
\Big\}
\]
d'un certain polynôme $R = R(x, y)$ de degré $d \geqslant \kappa+3$, 
pour tout ordre de jets intermédiaire:
\[
1
\,\leqslant\,
\lambda
\,\leqslant\,
\kappa,
\]
l'expression:
\[
\!\!\!\!\!\!\!\!\!\!\!\!\!\!\!\!\!\!\!\!
\!\!\!\!\!\!\!\!\!\!\!\!\!\!\!\!\!\!\!\!
\boxed{\,
{\sf J}_R^\lambda
\,:=\,
\left\{
\aligned
&
\ \ \ \ \,
\frac{y^{(\lambda)}}{R_x}
+\!\!
\sum_{\mu_1+\cdots+(\lambda-1)\mu_{\lambda-1}=\lambda}
\!\!\!\!\!\!\!
\frac{\big(y'\big)^{\mu_1}\cdots
\big(y^{(\lambda-1)}\big)^{\mu_{\lambda-1}}}{
R_x}\,
\mathcal{J}_{\mu_1,\dots,\mu_{\lambda-1}}^\lambda\!
\left(
\frac{R_y}{R_x},
\bigg(
\frac{R_{x^i y^j}}{R_x}
\bigg)_{2\leqslant i+j\leqslant
\atop
\leqslant
-1+\mu_1+\cdots+\mu_{\lambda-1}}
\right),
\\
\!
&
-\,
\frac{x^{(\lambda)}}{R_y}
-\!\!
\sum_{\mu_1+\cdots+(\lambda-1)\mu_{\lambda-1}=\lambda}
\!\!\!\!\!\!\!
\frac{\big(x'\big)^{\mu_1}\cdots\big(x^{(\lambda-1)}\big)^{\mu_{\lambda-1}}}{
R_y}\,
\mathcal{J}_{\mu_1,\dots,\mu_{\lambda-1}}^\lambda\!
\left(
\frac{R_x}{R_y},
\bigg(
\frac{R_{y^ix^j}}{R_y}
\bigg)_{2\leqslant i+j\leqslant
\atop
\leqslant
-1+\mu_1+\cdots+\mu_{\lambda-1}}
\right),
\\
&
\ \ \ \ \ \ \ \ \ \ \ 
0
\ \ \ \ \ \ \ \ \ 
\text{\rm sur}\ \
X^1\cap\P_\infty^1,
\endaligned
\right.}
\]
définit une \underline{\sl différentielle génératrice} de $\lambda$-jets
{\em holomorphe} sur $X^1$ tout entier qui s'annule identiquement sur
le diviseur ample $X^1 \cap \P_\infty^1$.\qed
\end{Theoreme}

\section{\bf Amplitude génératrice}
\label{amplitude-generatrice}
\HEAD{\ref{amplitude-generatrice}.~Amplitude génératrice}{
Jo\"el Merker, Département de Mathématiques d'Orsay}

La fin de l'énoncé du théorème donne des éléments qui sont
maintenant suffisants pour conclure, les vérifications
étant laissées au lecteur.
En particulier, les asymptotiques reposent sur des calculs
élémentaires, déjà connus (\cite{ Green-Griffiths-1980,
Merker-2010}), puisqu'en effet, il s'agit d'estimer:
\[
\!\!\!\!\!\!\!\!\!\!\!\!\!\!\!\!\!\!\!\!\!\!\!\!
\!\!\!\!\!\!\!\!\!\!\!\!\!\!\!\!\!\!\!\!\!\!\!\!
\aligned
\sum_{m_1+\cdots+\kappa m_\kappa=m}
\bigg\{
\binom{m_1(d-3)+\cdots+m_\kappa(d-\kappa-2)+2}{2}
-
\binom{m_1(d-3)+\cdots+m_\kappa(d-\kappa-2)-d+2}{2}
\bigg\},
\endaligned
\]
à savoir:
\[
\aligned
&
\equiv
\frac{2\,d}{2}
\!\!\!
\sum_{m_1+\cdots+\kappa m_\kappa=m}
\!\!\!
\Big\{
m_1(d-3)
+\cdots+
m_\kappa(d-\kappa-2)
\Big\}
+
m^\kappa\,{\rm O}(d)
+
{\rm O}\big(m^{\kappa-1}\big)
\\
&
\equiv
d\,
\sum_{m_1+\cdots+\kappa m_\kappa=m}
\!\!\!
\big(
m_1+\cdots+m_\kappa
\big)\,d
+
m^\kappa\,{\rm O}(d)
+
{\rm O}\big(m^{\kappa-1}\big),
\endaligned
\]
sachant que par approximation intégrale:
\[
\aligned
\sum_{m_1+\cdots+\kappa m_\kappa=m}
&
\equiv
\frac{m^\kappa}{\kappa!}\,
\int_{y_2+\cdots+y_\kappa\geqslant 1
\atop
y_1\geqslant 0,\,\dots,\,y_\kappa\geqslant 0}\,
\bigg(
1
-
\frac{y_2}{2}
-\cdots-
(\kappa-1)\,
\frac{y_\kappa}{\kappa}
\bigg)
+
{\rm O}\big(m^{\kappa-1}\big)
\\
&
=
\frac{m^\kappa}{\kappa!}\,
\bigg[
\frac{1}{(\kappa-1)!}
-
\frac{1}{2}\,\frac{1}{\kappa!}
-\cdots-
\frac{\kappa-1}{\kappa}\,
\frac{1}{\kappa!}
\bigg]
+
{\rm O}\big(m^{\kappa-1}\big)
\\
&
=
\frac{m^\kappa}{\kappa!\,\kappa!}\,
\bigg[
1
+
\frac{1}{2}
+\cdots+
\frac{1}{\kappa}
\bigg]
+
{\rm O}\big(m^{\kappa-1}\big),
\endaligned
\]
ce qui fournit bien le résultat asymptotique annoncé:
\[
\frac{m^\kappa}{\kappa!\,\kappa!}\,
\Big[
d^2\,\log\,\kappa
+
d^2\,{\rm O}(1)
+
{\rm O}(d)
\Big]
+
{\rm O}\big(m^{\kappa-1}\big).
\]

\medskip

En changeant de carte affine, la prescription d'annulation
à l'infini tourne. 

\medskip

Voir directement que le fibré $E_{\kappa, m}^{\rm GG} T_X^*$
est engendré par ses sections holomorphes globales
est alors aisé.

\medskip

Enfin, tout se généralise directement au cas des courbes
algébriques complexes géométriquement lisses:
\[
X^1
\,\subset\,
\P^{1+c}(\C)
\]
de codimension $c \geqslant 1$ quelconque qui sont
intersections complètes entre $c$ hypersurfaces
algébriques:
\[
\aligned
0
&
=
R^1\big(x_1,\dots,x_c,x_{c+1}\big),
\\
\cdot
&
\cdots\cdots\cdots\cdots\cdots\cdots\cdots\cdot
\\
0
&
=
R^c\big(x_1,\dots,x_c,x_{c+1}\big),
\endaligned
\] 
puisqu'alors les quantités $R_x$ et $R_y$ sont naturellement
remplacées par des mineurs généraux qui occupent des places
dénominatoriales:
\[
\frac{x_1'}{
\left\vert\!
\begin{array}{ccc}
R_{x_2}^1 & \cdots & R_{x_{c+1}}^1
\\
\cdot\cdot & \cdots & \cdot\cdot 
\\
R_{x_2}^c & \cdots & R_{x_{c+1}}^c
\end{array}
\!\right\vert}
\,=\,\cdots\cdots\,=\,
(-1)^c\,
\frac{x_{c+1}'}{
\left\vert\!
\begin{array}{ccc}
R_{x_1}^1 & \cdots & R_{x_c}^1
\\
\cdot\cdot & \cdots & \cdot\cdot 
\\
R_{x_1}^c & \cdots & R_{x_c}^c
\end{array}
\!\right\vert},
\]
les différentielles de jets génératrices d'ordre $\kappa \geqslant 1$
quelconque étant obtenues par différentiation directe de ces $c$
égalités.

\vfill\end{document}

%% file: macros.tex
\newcommand{\C}{\mathbb{C}}
\newcommand{\D}{\mathbb{D}}

\newcommand{\N}{\mathbb{N}}\renewcommand{\P}{\mathbb{P}}

\newcommand{\Z}{\mathbb{Z}}

\newcommand{\red}{\textcolor{red}}

\newcommand{\HEAD}[2]{%
\pagestyle{fancy}
\fancyhead[RO]{\tiny\sf\thepage}
\fancyhead[CO]{{\tiny\sf #1}}
\fancyhead[LE]{\tiny\sf\thepage}
\fancyhead[CE]{{\tiny\sf #2}}
\fancyfoot{}}

\newtheorem{The}{Theorem}[section]

\newtheorem{Theoreme}{Théorème}[section]
\newtheorem{Proposition}[The]{Proposition}
\newtheorem{Lemme}[The]{Lemme}

\theoremstyle{definition}

\newtheorem{Question}[The]{Question}

\newcommand{\Card}{\text{\scriptsize\sf Card}}

\renewcommand{\deg}{\text{\footnotesize\sf deg}}

\renewcommand{\dim}{\text{\footnotesize\sf dim}}

\newcommand{\infinibonus}{\text{\footnotesize\sf $\infty$-Bonus}}

\renewcommand{\lim}{\text{\footnotesize\sf lim}}

\renewcommand{\log}{\text{\footnotesize\sf log}}

\newcommand{\mathmotsf}[1]{\text{\footnotesize\sf #1}}

\newcommand{\ordreinfini}{\text{\footnotesize\sf Ordre-$\infty$}}

\newcommand{\Sym}{\text{\scriptsize\sf Sym}}

\newcommand{\vf}{\vfill\end{document}}

\newcommand{\zero}[1]{\underline{#1}_{\red{\circ}}}

\newcommand{\zerozero}[1]{\underline{#1}_{\red{\circ\circ}}}
